\documentclass[12pt]{amsart}
\RequirePackage{amssymb}
\RequirePackage{times}
\RequirePackage{graphicx}
\RequirePackage{hyperref}
\RequirePackage[hang,small,bf]{caption}

\newcommand{\R}{\mathbb{R}}

\newcommand{\E}{\mathbb{E}}
\newcommand{\Cov}{\mathrm{Cov}}

\newcommand{\Var}{\mathbb{V}\mathrm{ar}}

\newcommand\stack[2]{\genfrac{[}{]}{0pt}{}{#1}{#2}}

\def\Tr{\mathrm{Tr}}

\newtheorem{theorem}{Theorem}[section]
\newtheorem{lemma}[theorem]{Lemma}

\newtheorem{corollary}[theorem]{Corollary}

\makeatletter
\g@addto@macro{\endabstract}{\@setabstract}
\newcommand{\authorfootnotes}{\renewcommand\thefootnote{\@fnsymbol\c@footnote}}%
\makeatother

\makeatletter
\let\@fnsymbol\@arabic
\makeatother

\begin{document}
\begin{center}
  \LARGE 
 Limit theorems for linear eigenvalue statistics of overlapping matrices \par \bigskip
  \normalsize
  \authorfootnotes
  Vladislav Kargin\footnote{email: vladislav.kargin@gmail.com; current address: 282 Mosher Way, Palo Alto, CA 94304, USA} \par \bigskip
\end{center}

\begin{center}
\textbf{Abstract}
\end{center}

\begin{quotation}
The paper proves several limit theorems for linear eigenvalue statistics of
overlapping Wigner and sample covariance matrices. It is shown that the covariance of the limiting multivariate Gaussian distribution is
diagonalized by choosing the Chebyshev polynomials of the first kind as the
basis for the test function space. The covariance of linear statistics for the Chebyshev polynomials of sufficiently high degree depends only on the first two moments of the matrix entries. Proofs are based on a graph-theoretic interpretation of the Chebyshev linear statistics as sums over non-backtracking cyclic paths.
\end{quotation}

\section{Introduction}

Recently, there was a surge of interest in the spectral properties of
overlapping submatrices of large random matrices. A seminal study was done
by Baryshnikov \cite{baryshnikov01}, which derived the joint eigenvalue
distribution for principal minors of Gaussian Unitary Ensemble (GUE) and
related this distribution to the last passage percolation problem. Later
Johansson and Nordenstam \cite{johansson_nordenstam06} established the
determinantal structure of this joint distribution, and their results were
generalized in \cite{forrester_nagao08}, \cite{forrester_nordenstam09} and 
\cite{metcalfe11} to other unitarily invariant ensembles of random matrices.
Recently, Borodin in \cite{borodin10} and Reed in \cite{reed14} obtained
limit theorems for eigenvalue statistics of overlapping real Gaussian
matrices. We extend these results further to Wigner and sample covariance
matrices which lack the rotation invariance structure.

\bigskip

\subsection{Main Results}
\subsubsection{Wigner matrices}

\label{subsection_wigner_matrices}

\begin{figure}[btph]
\includegraphics[width=4.5cm]{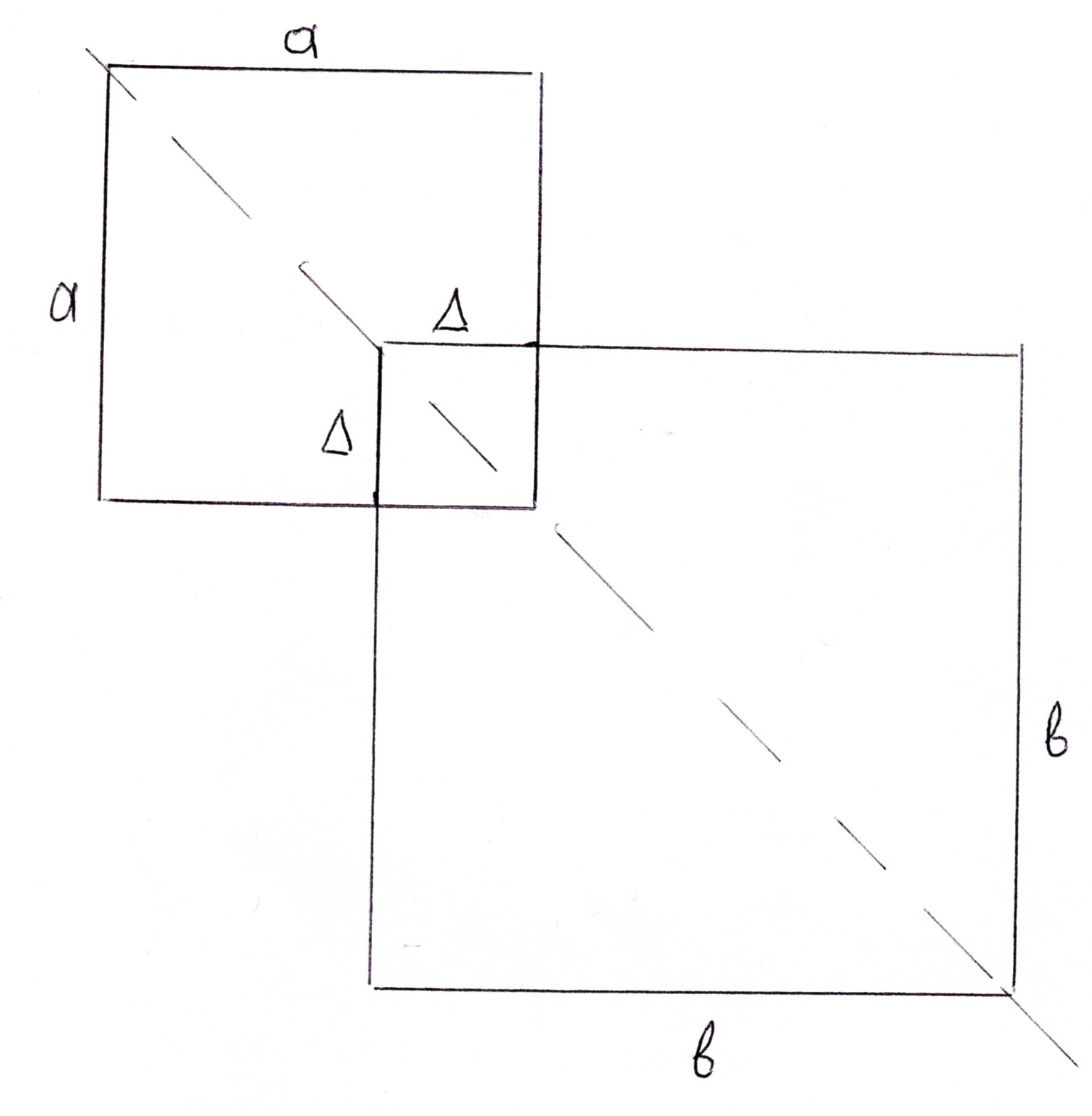} 
\caption{Overlapping symmetric matrices}
\label{fig:Wigner}
\end{figure}

The overlapping real Wigner matrices, $A_{N}$ and $B_{N}$, are two symmetric matrices
which are submatrices of an infinite real symmetric matrix $X.$ This matrix $X$
has independent identically distributed entries above the main diagonal and
independent entries (with a possibly different distribution) on the main
diagonal. The off-diagonal entries of $X$ are random variables with all moments
finite and $\mathbb{E}\left( X_{ij}\right) =0,$ $\mathbb{E}%
\left( X_{ij}^{2}\right) =1,$ and $\mathbb{E}\left( X_{ij}^{4}\right)
=m_{4}. $ The diagonal terms have all finite
moments and \ $\mathbb{E}\left( X_{ii}\right) =0$, $\mathbb{E}\left( X_{ii}^{2}\right) =d_{2}.$

Let $A_{N}$ have $a^{\left( N\right) }$ rows and columns, and $B_{N}$
have $b^{\left( N\right) }$ rows and columns. In addition, $A_{N}$ and $%
B_{N} $ have $\Delta ^{\left( N\right) }$ rows and columns in common. (See
Figure \ref{fig:Wigner}.)

We think about $A_{N}$ and $B_{N}$ as sequences of matrices of increasing
size and we assume that there is a parameter $t_{N}$ that approaches
infinity as $N\rightarrow \infty ,$ \ and that quantities $a^{(N)}/t_{N},$ $%
b^{\left( N\right) }/t_{N},$ and $\Delta ^{\left( N\right) }/t_{N}$ approach
some positive limits $a,$ $b,$ and $\Delta ,$ respectively.

 We define the normalized matrices 
\begin{equation*}
\widetilde{A}_{N}:=\frac{1}{2\sqrt{a^{\left( N\right) }}}A_{N},\text{ and }%
\widetilde{B}_{N}:=\frac{1}{2\sqrt{b^{\left( N\right) }}}B_{N}.
\end{equation*}%
The normalization is chosen in such a way that the empirical distribution of
the eigenvalues of $\widetilde{A}_{N}$ and $\widetilde{B}_{N}$ converges to a
distribution supported on the interval $\left[ -1,1\right] .$

If $f:\mathbb{R}\rightarrow \mathbb{R}$ is a test function, then 
the corresponding linear statistic for matrix $A_{N}$ is 
\begin{equation*}
\mathcal{N}\left( f,A_{N}\right) :=\mathrm{Tr}\left[ f\left( \widetilde{A}%
_{N}\right) \right] =\sum f\left( \lambda _{i}\right) ,
\end{equation*}%
where $\lambda _{i}$ are the eigenvalues of the matrix $\widetilde{A}_{N}.$
It is the linear statistic of the eigenvalues of $A_{N}$ after rescaling. We
define linear statistics for matrix $B_{N}$ similarly. Let us also define
the centered linear statistics, 
\begin{equation*}
\mathcal{N}^{o}\left( f,A_{N}\right) :=\mathcal{N}\left( f,A_{N}\right) -%
\mathbb{E}\mathcal{N}\left( f,A_{N}\right),
\end{equation*}
and similarly for $\mathcal{N}^{o}\left( f,B_{N}\right)$.

We are interested in the joint distribution of these linear statistics when $%
N$ is large. Recall that the Chebyshev polynomials of the first kind are
defined by the formula: $T_{k}\left( \cos \theta \right) =\cos k\theta .$ We
will prove the following result.

\begin{theorem}
\label{theorem_lin_statistics_Wigner_general} Let  $A_{N}$ and $B_{N}$ be overlapping real Wigner matrices. Let $T_{k}\left( x\right) $ denote
the Chebyshev polynomials of the first kind. As $N\rightarrow \infty ,$ the
distribution of the centered linear statistics $\mathcal{N}^{o}\left(
T_{k},A_{N}\right) $ and $\mathcal{N}^{o}\left( T_{l},B_{N}\right) $
converges to a two-variate Gaussian distribution with the covariance equal
to 
\begin{equation}
\left\{ 
\begin{array}{cc}
\frac{k}{2}\left( \frac{\Delta }{\sqrt{ab}}\right) ^{k}, & \text{if }k=l\geq
3, \\ 
\frac{m_{4}-1}{2}\left( \frac{\Delta }{\sqrt{ab}}\right) ^{2}, & \text{if }%
k=l=2, \\ 
\frac{d_{2}}{4}\left( \frac{\Delta }{\sqrt{ab}}\right) , & \text{if }k=l=1,
\\ 
0, & \text{otherwise.}%
\end{array}%
\right.  \label{formula_covariance1}
\end{equation}
\end{theorem}

For the author, the motivation for this model came from the paper by Borodin 
\cite{borodin10} where a similar problem was considered and solved for
matrices whose entries have the same first four moments as the Gaussian
random variable.  (The covariances of linear statistics depend only on the first four moments of matrix entries, as was noted by Anderson and Zeitouni in \cite{anderson_zeitouni06} and by Tao and Vu in 
\cite{tao_vu11}). 

Interestingly, Theorem \ref{theorem_lin_statistics_Wigner_general} shows that the fourth-order moment influences only the covariance of the Chebyshev polynomials of the second order. (For non-overlapping matrices, this fact can also be deduced from the formulas in Theorem 3.6 in
 \cite{lytova_pastur09} or Theorem 1 in \cite{shcherbina11}.) Our combinatorial proof sheds some light on the origin of this phenomenon.  

\subsubsection{Sample covariance matrices}

\begin{figure}[btph]
\includegraphics[width=4.5cm]{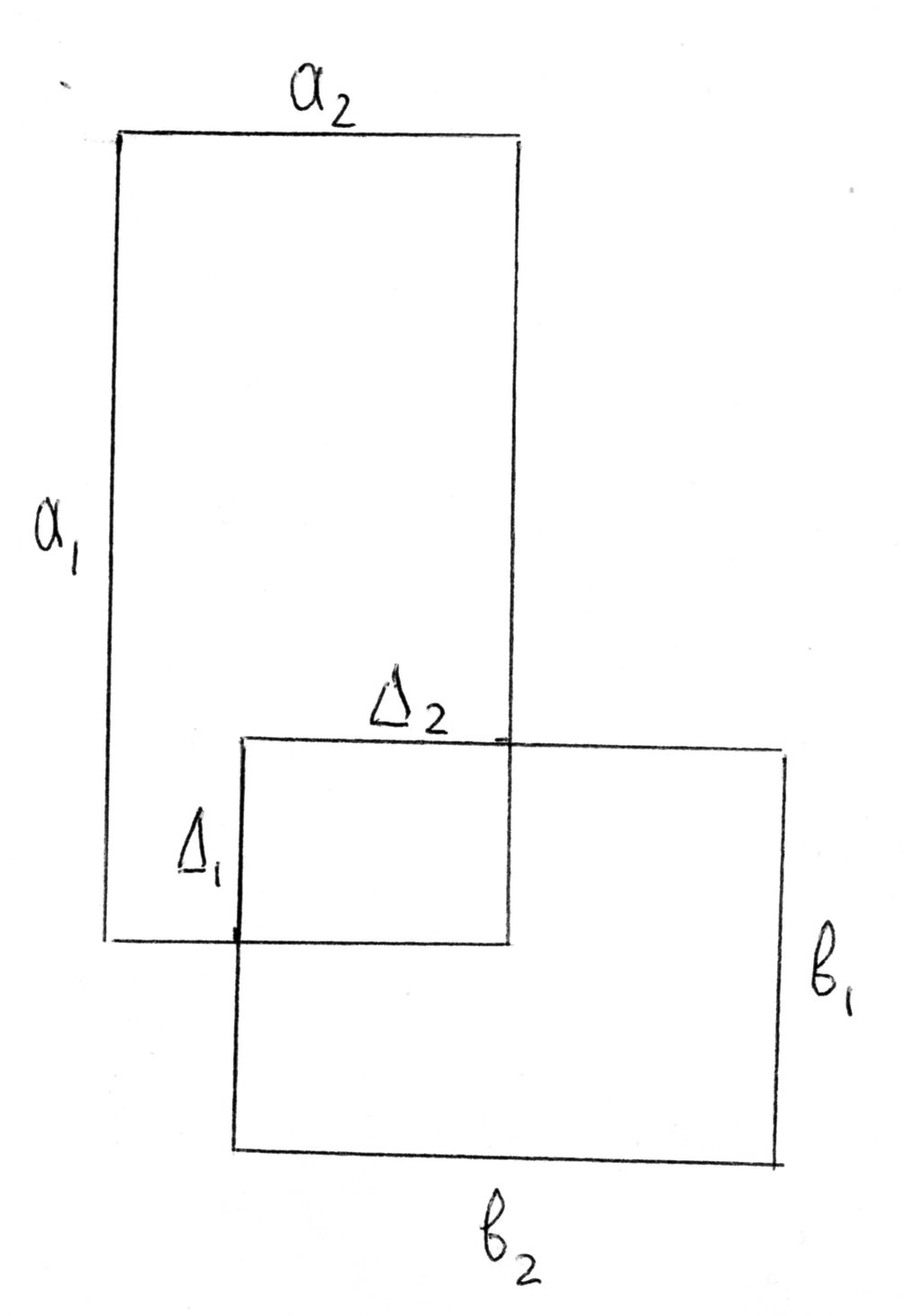} 
\caption{Overlapping non-symmetric matrices}
\label{fig:Wishart}
\end{figure}

Next, we study the singular values of non-symmetric overlapping matrices.
This is equivalent to the study of eigenvalues of certain sample covariance
matrices. Namely, let $A_{N}$ and $B_{N}$ be two submatrices of an infinite
random matrix $X$ with independent identically distributed random entries.
Let the entries of $X$ be real random variables with all moments finite and
assume that 
\begin{equation}
\E( X_{ij}) =0, \ \E(X_{ij}^2) =1, \text{ and }\E ( X_{ij}^4) =m_4.
\label{moments_Wishart}
\end{equation}
Suppose that $A_{N}$ has $a_{1}^{\left( N\right) }$ rows and $a_{2}^{\left(
N\right) }$ columns, and that $B_{N}$ has $b_{1}^{\left( N\right) }$ rows
and $b_{2}^{\left( N\right) }$ columns. In addition, suppose that $A_{N}$
and $B_{N}$ have $\Delta _{1}^{\left( N\right) }$ rows and $\Delta
_{2}^{\left( N\right) }$ columns in common. (See Figure \ref{fig:Wishart}.)

As before, we think about $A_{N}$ and $B_{N}$ as sequences of matrices of
increasing size and we assume that there is a parameter $t_{N}$ that
approaches infinity as $N\rightarrow \infty ,$ \ and that quantities $%
a_{1}^{(N)}/t_{N},$ $a_{2}^{\left( N\right) }/t_{N},$ $b_{1}^{\left(
N\right) }/t_{N},$ $b_{2}^{\left( N\right) }/t_{N},$ $\Delta _{1}^{\left(
N\right) }/t_{N},$ and $\Delta _{2}^{\left( N\right) }/t_{N}$ approach some
positive limits $a_{1},$ $a_{2},$ $b_{1},$ $b_{2},$ $\Delta _{1},$ and $%
\Delta _{2},$ respectively.

Let us define the normalized sample covariance matrices 
\begin{equation}
W_{A_{N}}:=\frac{1}{2\sqrt{a_{1}^{\left( N\right) }a_{2}^{(N)}}}\left[
A_{N}A_{N}^{\ast }-\left( a_{1}^{\left( N\right) }+a_{2}^{\left( N\right)
}\right) I_{a_{1}^{\left( N\right) }}\right] , \label{W_A_definition}
\end{equation}%
where $I_{a_{1}^{\left( N\right) }}$ is the $a_{1}^{\left( N\right) }$-by-$%
a_{1}^{\left( N\right) }$ identity matrix, and 
\begin{equation}
W_{B_{N}}:=\frac{1}{2\sqrt{b_{1}^{\left( N\right) }b_{2}^{(N)}}}\left[
B_{N}B_{N}^{\ast }-\left( b_{1}^{\left( N\right) }+b_{2}^{\left( N\right)
}\right) I_{b_{1}^{\left( N\right) }}\right] . \label{W_B_definition}
\end{equation}%
Again, the normalization is chosen in such a way that the empirical
distribution of the eigenvalues of $W_{A_{N}}$ and $W_{B_{N}}$ converges to a
distribution supported on the interval $\left[ -1,1\right] .$

If $f:\mathbb{R}\rightarrow \mathbb{R}$ is a test function, then we define
the corresponding linear statistic for matrix $W_{A_{N}}$ as 
\begin{equation*}
\mathcal{N}\left( f,W_{A_{N}}\right) :=\mathrm{Tr}\left[ f\left(
W_{A_{N}}\right) \right] =\sum f\left( \lambda _{i}\right) ,
\end{equation*}%
where $\lambda _{i}$ are the eigenvalues of the matrix $W_{A_{N}}.$ Hence,
this is a linear statistic of the matrix $A_{N}$'s squared singular
values. We define linear statistics for matrix $W_{B_{N}}$ similarly. Let us
also define the centered linear statistics, 
\begin{equation}
\mathcal{N}^{o}\left( f,W_{A_{N}}\right) :=\mathcal{N}\left(
f,W_{A_{N}}\right) -\mathbb{E}\mathcal{N}\left( f,W_{A_{N}}\right) ,  \label{cent_stat_Wishart} 
\end{equation}%
and similarly for $\mathcal{N}^{o}\left( f,W_{B_{N}}\right)$. Finally, let 
\begin{equation*}
\gamma :=\frac{\Delta _{1}\Delta _{2}}{\sqrt{a_{1}a_{2}b_{1}b_{2}}}.
\end{equation*}

\begin{theorem}
\label{theorem_lin_statistics_Wishart_general} Assume $A_{N}$ and $B_{N}$
are the real random overlapping matrices with the matrix entries that satisfy (\ref{moments_Wishart}).
 Let $T_{k}\left( x\right) $ denote the
Chebyshev polynomials of the first kind. As $N\rightarrow \infty ,$ the
distribution of the centered linear statistics $\mathcal{N}^{o}\left(
T_{k},W_{A_{N}}\right) $ and $\mathcal{N}^{o}\left( T_{l},W_{B_{N}}\right) $
converges to a two-variate Gaussian distribution with the covariance equal
to 
\begin{equation*}
\left\{ 
\begin{array}{cc}
\delta _{kl}\frac{k}{2}\gamma ^{k}, & \text{if }k>1, \\ 
\delta _{kl}\frac{\left( m_{4}-1\right) }{4}\gamma , & \text{if }k=1.%
\end{array}%
\right.
\end{equation*}
\end{theorem}

This result parallels the result for Wigner matrices and has the same surprising conclusion that the fourth moment of matrix entries influences only the covariances of low-degree Chebyshev polynomials. In this case, these are the first order polynomials.

\subsubsection{CLT for continuously differentiable functions}

In order to extend our results to continuously differentiable functions, we
have to restrict to models with matrix entries that satisfy
the Poincare inequality property. 
 
 We say that a matrix entry $X_{ij}$ satisfies the Poincare inequality property if there is a
constant $c>0$ such that for every continuously differentiable function $%
f\left( x\right) ,$ we have 
\begin{equation}
\mathbb{V}ar\left( f\left( X_{ij}\right) \right) \leq c\mathbb{E}\left( |
\nabla f\left( X_{ij}\right) |^{2}\right) .  \label{Poincare_inequality}
\end{equation}
For example, the Poincare inequality property holds for models with Gaussian entries 
or with entries uniformly distributed on the unit interval but not for
the model with $\pm 1$ entries. 

We will use the Poincare inequality to bound variances of linear statistics of large dimensional matrices and prove the tightness of their distributions.  We consider only the case of
Wigner matrices. The results can be extended to the case of sample covariance matrices.

First, let us define the coefficients in the expansion of a function $f$
over Chebyshev polynomials: 
\begin{equation}
\widehat{f}_{k}:=\left\{ 
\begin{array}{cc}
\frac{2}{\pi }\int_{-1}^{1}f\left( x\right) T_{k}\left( x\right) \frac{dx}{%
\sqrt{1-x^{2}}}, & \text{for }k\geq 1, \\ 
\frac{1}{\pi }\int_{-1}^{1}f\left( x\right) \frac{dx}{\sqrt{1-x^{2}}}, & 
\text{for }k=1.%
\end{array}%
\right.  \label{Chebyshev_coefficients}
\end{equation}

Let $\mathcal{F}$ be the linear subspace of functions $f:\mathbb{R}%
\rightarrow \mathbb{R}$ which are differentiable with continuous derivative
in the interval $I_{\delta }=\left[ 1-\delta ,1+\delta \right],$ and grow no faster than a polynomial outside of this interval.

\begin{theorem}
\label{theorem_lin_statistics_diff_functions}Assume  $A_{N}$ and $B_{N}$
are overlapping real Wigner matrices, and their matrix entries satisfy the Poincare inequality. Then for every pair of
functions $f$ and $g$ in $\mathcal{F}$, the linear statistics $\mathcal{N}%
^{o}\left( f,A_{N}\right) $ and $\mathcal{N}^{o}\left( g,B_{N}\right) $
converge in distribution to the bivariate Gaussian random variable with the
covariance matrix $V:$ 
\begin{eqnarray*}
V_{11} &=&\frac{1}{2}\left( \frac{d_{2}}{2}\left( \widehat{f}_{1}\right)
^{2}+\left( m_{4}-1\right) \left( \widehat{f}_{2}\right)
^{2}+\sum_{k=3}^{m}k\left( \widehat{f}_{k}\right) ^{2}\right) \\
V_{12} &=&\frac{1}{2}\left( \frac{d_{2}}{2}\widehat{f}_{1}\widehat{g}%
_{1}\gamma +\left( m_{4}-1\right) \widehat{f}_{2}\widehat{g}_{2}\gamma
^{2}+\sum_{k=3}^{m}k\widehat{f}_{k}\widehat{g}_{k}\left( \gamma \right)
^{k}\right) \\
V_{22} &=&\frac{1}{2}\left( \frac{d_{2}}{2}\left( \widehat{g}_{1}\right)
^{2}+\left( m_{4}-1\right) \left( \widehat{g}_{2}\right)
^{2}+\sum_{k=3}^{m}k\left( \widehat{g}_{k}\right) ^{2}\right) ,
\end{eqnarray*}%
where $\gamma =\Delta /\sqrt{ab}.$
\end{theorem}

\subsection{Discussion}
Let us put our results in a more general prospective. Linear
statistics of sample covariance matrices have been first investigated in
Jonsson \cite{jonsson82}, who established the joint CLT for the moments of the empirical eigenvalue distribution. A recent contribution by Anderson and Zeitouni \cite%
{anderson_zeitouni06} extended the study of linear statistics to a very
general class of matrices with independent entries of changing variance.
They derived a formula for the covariance of linear eigenvalue statistics and
proved a CLT theorem for continuously differentiable test functions when the matrix entries satisfy the Poincare inequality. They have also
noted a relation to the Chebyshev polynomials. Their method is combinatorial, in the spirit of the method of moments. 

For more restricted classes
of matrices, namely, for Gaussian and unitarily invariant matrices,
important results about linear statistics were established in
Diaconis-Shahshahani \cite{diaconis_shahshahani94}, Costin-Lebowitz \cite%
{costin_lebowitz95}, Johansson \cite{johansson98}, Soshnikov \cite%
{soshnikov00} and \cite{soshnikov00a}, and Diaconis-Evans \cite%
{diaconis_evans01}. 

In \cite{lytova_pastur09}, Lytova and Pastur 
showed how to use an analytic method to interpolate the results from Gaussian and Wishart ensembles
to Wigner and sample-covariance matrices. For example, they derived a CLT theorem for test functions in the class $\mathcal{C}^5$
 (five continuously differentiable derivatives) when the matrix entries have 5 finite moments. Later, the proofs were simplified and conditions on the test function smoothness were weakened in \cite{shcherbina11}.

The main novelty of our results is that they extend the investigation to the
spectra of overlapping submatrices. This was started by Borodin in \cite{borodin10} who investigated the Gaussian case of Wigner matrices. Here we treat the case of non-Gaussian Wigner and sample covariance matrices.

 The second contribution is that we give a combinatorial explanation for the important
role of the Chebyshev polynomials of the first kind. These polynomials were related 
to the fluctuations of the empirical eigenvalue distribution in Johansson 
 (Corollary 2.8 in \cite{johansson98}) by analytical methods, and this relation was later extended to a dynamical version in Cabanal-Duvillard \cite{cabanal_duvillard01}. 
Feldheim and Sodin in \cite{feldheim_sodin10} explained the combinatorial significance of the Chebyshev polynomials of
the second type in the context of random matrices, and we extend their work by explaining the role of the Chebyshev polynomials of the first type.

Finally, recall
that from the \textquotedblleft four moment theorem\textquotedblright\ by
Tao and Vu \cite{tao_vu11}, we can expect that the CLT covariance matrix for linear statistics depends only on the first four moments of the
matrix entries. For non-overlapping matrices,  Anderson and Zeitouni in \cite{anderson_zeitouni06} and Lytova and Pastur  in  \cite{lytova_pastur09} derived explicit formulas for the CLT covariance matrix by using combinatorial and analytic methods, respectively. In particular, they showed that the third moment of matrix entries does not influence the covariance of linear statistics. We extend these findings to overlapping matrices and find an additional interesting fact that in the basis of the Chebyshev $T$-polynomials, the fourth moment affects only the covariance of the second or first degree polynomials, for the Wigner and sample
covariance matrices, respectively. 

The rest of the paper is organized as follows. In Section \ref{section_NBT},
we will show how the Chebyshev polynomials of the first type are related to
the non-backtracking tailless paths on regular and bi-regular graphs. In Section \ref{section_polynomials},
we will prove the CLT results for simple models in which the entries are
either $\pm 1$ or uniformly distributed on the unit circle (Theorems \ref%
{theorem_lin_statistics_Wigner} and \ref{theorem_lin_statistics_Wishart}).
In Section \ref{section_general_matrices} we will prove Theorems \ref%
{theorem_lin_statistics_Wigner_general} and \ref%
{theorem_lin_statistics_Wishart_general}. Section \ref{section_C1} is
devoted to the proof of Theorem \ref{theorem_lin_statistics_diff_functions}.
And Section \ref{section_conclusion} concludes.

\section{The Chebyshev polynomials and non-backtracking paths}

\label{section_NBT}

The Chebyshev polynomials of the first and the second kind are defined by
the formulas $T_k( \cos \theta ) =\cos ( k\theta ) $
and $U_{k}( \cos \theta ) =\sin ((k+1)\theta) /\sin
\theta ,$ respectively. For $k\geq 2,$ both the $T$ and $U$ polynomials
satisfy the same recursion. For example, for $U$ polynomials, it is $%
U_{k}\left( x\right) =2xU_{k-1}\left( x\right) -U_{k-2}\left( x\right) .$
The inital conditions in this recursion are $T_{0}(x)=1$ and $T_{1}\left(
x\right) =x$ for $T$-polynomials and $U_{0}(x)=1$ and $U_{1}\left( x\right)
=2x$ for $U$-polynomials.

The $T$ and $U$ polynomials can be related as follows: 
\begin{equation}
U_{k}\left( x\right) =2\left[ T_{k}\left( x\right) +T_{k-2}\left( x\right)
+...+T_{\varepsilon }\left( x\right) \right] +\left( \varepsilon -1\right) ,
\label{relation_U_and_T}
\end{equation}%
where $\varepsilon =1$ if $k$ is odd and $\varepsilon =0$ if $k$ is even.

\subsection{Regular graphs and Chebyshev polynomials}

Let $G$ be a $\left( d+1\right) $-regular graph with the vertex set $%
V=\left\{ 1,\ldots ,n\right\} .$ We will say that the matrix $A$ is a \emph{%
generalized adjacency matrix} if it is Hermitian and if $\left\vert
A_{uv}\right\vert =1$ if $\left( u,v\right) $ is an edge of $G$ and $0$
otherwise.

A path $\gamma $ is a sequence of vertices $\left( u_{0},u_{1},\ldots
,u_{k}\right) $ which are adjacent in graph $G.$ The length of this path is $%
k.$ The path is called \emph{non-backtracking} if $u_{j+1}\neq u_{j-1}.$ It
is \emph{closed} if $u_{k}=u_{0}.$ A closed path is \emph{non-backtracking
tailless} if it is non-backtracking and $u_{k-1}\neq u_{1}.$

Theorem \ref{theorem_NBT_U_Wigner} and the following remark are essentially
due to Feldheim and Sodin \cite{feldheim_sodin10}. (See Lemma II.1.1 and
Claim II.1.2 in their paper, where this result is proved for complete
graphs.)

\begin{theorem}[Feldheim-Sodin]
\label{theorem_NBT_U_Wigner}Suppose that $A$ is a generalized adjacency
matrix for a $\left( d+1\right) $-regular graph $G.$ Then for all $k\geq 1,$ 
\begin{equation}
\sum_{u_{0}=u,u_{k}=v}A_{u_{0}u_{1}}A_{u_{1}u_{2}}\ldots A_{u_{k-1}u_{k}}= 
\left[ P_{k}\left( A\right) \right] _{uv}  \label{P_and_NB_paths}
\end{equation}%
where the sum is over all non-backtracking paths of length $k$ from $u$ to $v,$
and $P_{k}\left( x\right) $ is a polynomial $\ $defined for $k=1,2$ as $%
P_{1}\left( x\right) :=x$, $P_{2}\left( x\right) :=x^{2}-\left( d+1\right) $,
and for $k\geq 3$ by the recursion: 
\begin{equation}
P_{k}\left( x\right) :=xP_{k-1}\left( x\right) -dP_{k-2}\left( x\right) .
\label{recursion_P}
\end{equation}
\end{theorem}

\textbf{Remark:} $P_{k}\left( x\right) $ can be expressed in terms of
Chebyshev's $U$-polynomials as follows:%
\begin{equation}
P_{k}\left( x\right) =d^{k/2}U_{k}\left( \frac{x}{2\sqrt{d}}\right)
-d^{\left( k-2\right) /2}U_{k-2}\left( \frac{x}{2\sqrt{d}}\right) .
\label{relation_P_and_U}
\end{equation}

Let us use the following notation
\begin{equation}
A(\gamma):=A_{u_{0}u_{1}}A_{u_{1}u_{2}}\ldots A_{u_{k-1}u_{k}},
\label{def_A1}
\end{equation}
 where $\gamma$ is a path $( u_{0},u_{1},\ldots ,u_{k-1},u_{k})$.

\begin{theorem}
\label{theorem_NBT_T_Wigner}
Suppose that $A$ is a generalized adjacency
matrix for a $\left( d+1\right) $-regular graph $G.$ Then for all $k\geq 1,$%
\begin{equation*}
\sum A(\gamma)=\left\{ 
\begin{array}{cc}
\mathrm{Tr}\left[ 2d^{k/2}T_{k}\left( \frac{A}{2\sqrt{d}}\right) \right] , & 
\text{ if }k\text{ is odd,} \\ 
\mathrm{Tr}\left[ 2d^{k/2}T_{k}\left( \frac{A}{2\sqrt{d}}\right) +\left(
d-1\right) I_{n}\right] , & \text{if }k\text{ is even,}%
\end{array}%
\right.
\end{equation*}%
where the sum on the left hand-side is over all closed non-backtracking
tailless paths $\gamma$ of length 
$k.$
\end{theorem}

\textsc{Proof of theorem 2.2.} First we say that a closed non-backtracking path $(u_0,u_1, \ldots, u_{k-1},u_0)$ has a tail of
length $l$ whenever
\begin{equation*}
u_0 = u_k, \  u_1 = u_{k-1}, \ \ldots,   u_l = u_{k-l}, \text{ and} \ u_{l+1} \ne u_{k-l-1}.
\end{equation*}
A tailless path, therefore, has a tail of length 0.

Define,
\begin{equation*}
Q_k(A) := \sum A(\gamma),
\end{equation*}
where the sum is over all closed non-backtracking tailless paths of length $k$. Note, Theorem \ref{theorem_NBT_U_Wigner} gives,
\begin{equation*}
\Tr[P_k(A)] = \sum A(\gamma),
\end{equation*}
where the sum is over all closed non-backtracking paths of length $k$. We partition the
sum depending on the tail length to get,
\begin{equation*}
\Tr[P_k(A)] = \sum A(\gamma) + \sum A(\gamma) + \ldots 
\end{equation*}
where the first term on the r.h.s. is the sum over all closed non-backtracking paths with a tail of length 0, the second term is the sum over all closed non-backtracking paths with a tail of length 1, etc.

The first term on the r.h.s. is $Q_k(A)$. Consider the second term on the r.h.s. Recall that each $\gamma=(u_0,\ldots,u_k)$ is non-backtracking with a tail of length
 $1$. This is true if and only if $(u_1, \ldots, u_{k-1})$ is a closed non-backtracking tailless path of length $k-2$, $u_0 = u_k$, $u_0 \ne u_2$
and $u_k \ne u_{k-2}$. Thus, given a fixed $(u_1, \ldots, u_{k-1})$ which is closed non-backtracking tailless,
there are $d-1$ choices for the tail $(u_0, u_1) = (u_k, u_{k-1})$ 
($u_1 = u_{k-1}$ has $d + 1$ neighbors, but $u_0 \ne u_2$ and $u_k \ne u_{k-2}$). Thus the second term equals $(d - 1)Q_{k-2}(A)$. Similar
considerations show that the third term on the r.h.s. equals $d(d-1)Q_{k-4}(A)$, the fourth term equals $d^2(d-1)Q_{k-6}(A)$, the fifth term equals $d^3(d-1)Q_{k-8}(A)$, etc. Therefore,
\begin{equation*}
\Tr[P_k(A)] = Q_k(A) + (d-1)Q_{k-2}(A) + d(d - 1)Q_{k-4}(A) + d^2(d - 1)Q_{k-6}(A) + \ldots
\end{equation*}

Consider the l.h.s. Recall that
\begin{eqnarray}
P_1(x) & = & x = 2\sqrt{d}\ T_1\left(\frac{x}{2\sqrt{d}}\right ), \notag \\
P_2(x) & = & x^2 - (d+1) = 2d \ T_2\left(\frac{x}{2\sqrt{d}}\right ) + (d-1). \notag
\end{eqnarray}
Also, whenever $k \geq 3$ and odd, equations (\ref{relation_U_and_T}) and (\ref{relation_P_and_U}) give,
\begin{equation}
P_k(x) = 2d^{k/2}T_k\left(\frac{x}{2\sqrt{d}}\right )
+ d^{(k-2)/2}(d-1)\left[2T_{k-2}\left(\frac{x}{2\sqrt{d}}\right )
+ \cdots + 2T_1\left(\frac{x}{2\sqrt{d}}\right ) \right]. \notag
\end{equation}
Finally, whenever $k \geq 3$ and even, equations (\ref{relation_U_and_T}) and (\ref{relation_P_and_U}) give,
\begin{equation*}
P_k(x) = 2d^{k/2}T_k\left(\frac{x}{2\sqrt{d}}\right )
+ d^{(k-2)/2}(d-1)\left [ 2T_{k-2}\left(\frac{x}{2\sqrt{d}}\right )
+ \cdots + 2T_0\left(\frac{x}{2\sqrt{d}}\right ) - 1\right ].
\end{equation*}
A proof by induction then gives the required result. \hfill $\square$

\subsection{Bipartite biregular graphs and Chebyshev polynomials}

A bipartite graph is a graph whose vertices belong to two sets, $V$ and $W,$
such that the vertices in $V$ are connected only to vertices in $W$ and vice
versa. A bipartite graph is $\left( c+1,d+1\right) $-regular if every vertex
in $V$ is connected to $c+1$ vertices in $W,$ and every vertex in $W$ is
connected to $d+1$ vertices in $V.$

Let $G$ be a $\left( c+1,d+1\right) $-regular graph with $\left\vert
V\right\vert =n$ and $\left\vert W\right\vert =m.$ Consider an $n$-by-$m$
matrix $A.$ We will identify row indices with elements of $V$ and column
indices with elements of $W.$ We say that an $n$-by-$m$ matrix $A$ is a
generalized adjacency matrix for a bipartite graph $G,$ if $\left\vert
A_{uv}\right\vert =1$ for $\left( u,v\right) \in G$ and $A_{uv}=0$ otherwise.

We define
\begin{equation}
A_2(\gamma):=A_{v_{0}w_{1}}\overline{A_{v_{1}w_{1}}}\ldots A_{v_{k-1}w_{k}}\overline{A_{v_{k}w_{k}}},
\label{def_A2}
\end{equation}
 where $\gamma$ is a path $(v_{0},w_{1},v_{1},\ldots ,w_{k},v_{k}) $ .

Let us define the following quantity%
\begin{equation}
\label{R_definition}
R_{k}\left( A,v_{0},v_{k}\right)
 =\sum A_2(\gamma),  
\end{equation}%
where the summation is over all non-backtracking paths $\gamma$ of length $2k$ from $v_{0}$
to $v_{k}.$

The following two results are essentially due to Feldheim and Sodin \cite%
{feldheim_sodin10}. (However, their expressions for $R_{k}$ in terms of
Chebyshev polynomials are incorrect. Compare their Lemma IV.1.1 and Claim
IV.1.2.)

\begin{theorem}
Suppose that the matrix $A$ is a generalized adjacency matrix for a $\left(
c+1,d+1\right) $-regular bipartite graph $G.$ Then for all $k\geq 1,$ 
\begin{equation}
R_{k}\left( A,v_{0},v_{k}\right) =\left[ F_{k}\left( AA^{\ast }\right) %
\right] _{v_{0}v_{k}},
\end{equation}%
where $F_k(x) $ are polynomials, which for $k=1,2$ \ are
defined as $F_1( x) :=x-( c+1) $ and $F_2(x) :=x^2-( 2c+d+1) x+( c+1) c,$ and for $k\geq
3 $, by the following recursion:%
\begin{equation*}
F_k( x) :=(x-( c+d) )F_{k-1}( x)
-cdF_{k-2}(x) .
\end{equation*}
\label{theorem_NBT_U_Wishart}
\end{theorem}

Note that if we define  \begin{equation*}
\widetilde{U}_{k}\left( x\right) :=\left( cd\right) ^{k/2}U_{k}\left( \frac{%
x-\left( c+d\right) }{2\sqrt{cd}}\right) ,
\end{equation*}%
where $U_{k}$ are the Chebyshev polynomials of the second kind and $U_{k}:=0$
for $k<0,$ then, for $k\geq 1$,
\begin{equation}
F_k(x) =\widetilde{U}_k(x) +( d-1) 
\widetilde{U}_{k-1}( x) -d\widetilde{U}_{k-2}( x) .
\label{relation_F_and_U}
\end{equation}
This can be checked by verifying the recursion for $F_k(x)$.

\textsc{Proof of Theorem 2.3.} It is easy to check the statement for $k = 1$. Indeed,
\begin{equation*}
[F_1(AA^{\ast})]_{v_0 v_1} = [AA^{\ast} - (c + 1)I]_{v_0 v_1} =
\left\{
\begin{array}{ll}
\sum_{w_1} A_{v_0 w_1} \overline{A_{v_1 w_1}}, & \text{if } v_0\ne v_1, \\
\sum_{w_1} A_{v_0 w_1} \overline{A_{v_1 w_1}} - (c+1), & \text{if } v_0 = v_1. 
\end{array}
\right .
\end{equation*}
Thus, since $A$ is an adjacency matrix of a $(c + 1, d + 1)$-regular bipartite graph,
\begin{equation*}
[F_1(AA^{\ast})]_{v_0 v_1} = \left\{
\begin{array}{ll}
\sum_{(v_0,w_1,v_1)} A_{v_0 w_1} \overline{A_{v_1 w_1}}, & \text{if } v_0\ne v_1, \\
0, & \text{if } v_0 = v_1, 
\end{array}
\right .
\end{equation*}
where the sum, when $v_0 \ne v_1$, is over all paths $(v_0,w_1, v_1)$. Note that, when $v_0 \ne v_1$, all
paths $(v_0, w_1, v_1)$ are necessarily non-backtracking. Also, when $v_0 = v_1$, there are no non-backtracking paths $(v_0,w_1, v_1)$. Therefore 
$[F_1(AA^{\ast})]_{v_0v_1} = R_1(A, v_0, v_1)$, by definition.

\begin{figure}[btph]
\includegraphics[width=15cm]{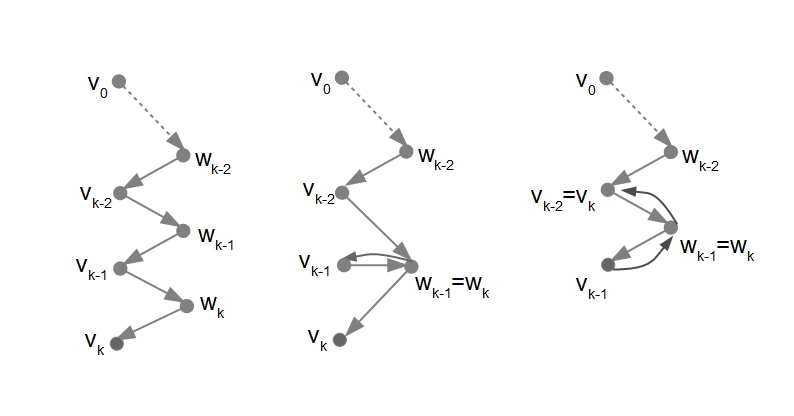} 
\caption{Three possibilities}
\label{fig:possibilities}
\end{figure}

For $k \geq 2$, consider
\begin{equation*}
\sum A_2(\gamma) ,
\end{equation*}
where the sum is over all paths $\gamma=(v_0,w_1,v_1,\ldots,w_k,v_k)$ for which $(v_0,w_1,v_1,\ldots,w_{k-1},v_{k-1})$
and $(v_{k-1},w_k, v_k)$ are both non-backtracking. There are three possibilities for such paths:
\begin{itemize}
\item $w_{k-1} \ne w_k$.
\item $w_{k-1} = w_k$ and $v_{k-2} \ne v_k$.
\item $w_{k-1} = w_k$ and $v_{k-2} = v_k$.
\end{itemize}
These possibilities are illustrated in Figure \ref{fig:possibilities}.

The first possibility is satisfied if $(v_0, w_1, v_1, \ldots, w_k, v_k)$ is non-backtracking. Thus the sum
over all terms which satisfy the first possibility is $R_k(A, v_0, v_k)$.

 The second possibility is
satisfied if  $(v_0,w_1, v_1, \ldots, w_{k-1}, v_k)$ is a non-backtracking path of length $2(k-1)$, $v_{k-2} \ne v_{k-1}$ and $v_{k-1} \ne v_k$. Thus, summing over all terms which satisfy the second possibility,
each fixed $(v_0, w_1, v_1, \ldots, w_{k-1}, v_k)$ which is non-backtracking is included $d-1$ times (there
are $d-1$ choices for the edge $(w_{k-1}, v_{k-1})$ since $w_{k-1}$ has 
$d + 1$ neighbors and $v_{k-1} \notin \{v_{k-2}, v_k\}$). Therefore the sum over all terms which satisfy the second possibility is 
$(d-1)R_{k-1}(A, v_0, v_k)$.

The third possibility is satisfied if $(v_0,w_1, v_1, \ldots, w_{k-2}, v_{k-2})$ is a non-backtracking
path of length $2(k - 2)$, $(v_{k-2}, w_{k-1}, v_{k-1})$ is a non-backtracking path, $w_{k-2} \ne w_{k-1}$, and
$(v_{k-1},w_k, v_k)$ is the path $(v_{k-2},w_{k-1}, v_{k-1})$ in reverse. Thus, summing over all terms which
satisfy the third possibility, each fixed $(v_0,w_1, v_1, \ldots, w_{k-2}, v_{k-2})$ which is non-backtracking
is included $(c + 1)d$ times when $k = 2$ and $cd$ times when $k \geq 3$. (These are the number of choices of the path $(v_{k-2},w_{k-1}, v_{k-1})$ with the above restrictions. To see this, note $v_{k-2}$ is fixed and has 
$c + 1$ neighbors. Thus, since $w_{k-2} \ne w_{k-1}$, there are $c + 1$ choices
for $w_{k-1}$ when $k = 2$ and $c$ choices for $w_{k-1}$ when $k \geq 3$. Also given such a choice
for $w_{k-1}$, note that $w_{k-1}$ has $d + 1$ neighbors. Thus, since $v_{k-1} \ne v_{k-2}$, there are $d$
choices for $v_{k-1}$.) Thus the sum over all the terms which satisfy the third possibility is
$(c+1)dR_{k-2}(A, v_0, v_{k-2}) = (c+1)dR_{k-2}(A, v_0, v_k)$ when 
$k = 2$, and $cdR_{k-2}(A, v_0, v_k)$ when
$k \geq 3$.

Alternatively note that, by definition,
\begin{equation*}
\sum A_2(\gamma) =
\sum_{v_{k-1}}
R_{k-1}(A, v_0, v_{k-1})R_1(A, v_{k-1}, v_k),
\end{equation*}
where the sum on the l.h.s. is over all paths for which $(v_0,w_1, v_1, \ldots,w_{k-1}, v_{k-1})$ and
$(v_{k-1},w_k, v_k)$ are both non-backtracking, and the sum on the r.h.s. is over all $v_{k-1} \in V$ .
For $k \geq 2$, the above observations thus give,
\begin{eqnarray}
\sum_{v_{k-1}}
R_{k-1}(A, v_0, v_{k-1})R_1(A, v_{k-1}, v_k) \notag \\
= R_k(A, v_0, v_k) + (d - 1)R_{k-1}(A, v_0, v_k) + (c + 1_{k=2})dR_{k-2}(A, v_0, v_k). \notag
\end{eqnarray}
A proof by induction then gives the required result. \hfill $\square$

Now, let us define 
\begin{equation}
\widetilde{T}_{k}\left( x\right) :=\left( cd\right) ^{k/2}T_{k}\left( \frac{%
x-\left( c+d\right) }{2\sqrt{cd}}\right) ,  \label{definition_T_modified}
\end{equation}%
where $T_{k}$ are the Chebyshev polynomials of the first kind$.$

\begin{theorem}
\label{theorem_NBT_T_Wishart}Suppose that the matrix $A$ is a generalized
adjacency matrix for a $\left( c+1,d+1\right) $-regular bipartite graph $G$
with vertex set $V\cup W.$ Assume $c\geq d.$ Then for all $k\geq 1,$%
\begin{equation*}
\sum A_2(\gamma)=\mathrm{Tr}%
\left[ 2\widetilde{T}_{k}\left( AA^{\ast }\right) +s_{c,d,k}I\right] 
\end{equation*}%
where the sum is over all closed non-backtracking tailless paths of length $2k$
that start with a vertex in $V$, and 
\begin{equation*}
s_{c,d,k}=\frac{( c-d)( -d)^k+cd-1}{d+1}.
\end{equation*}
\end{theorem}

Proof of Theorem \ref{theorem_NBT_T_Wishart} is similar to the proof of Theorem \ref{theorem_NBT_T_Wigner} and is deferred to Appendix.

\section{Linear statistics of polynomials for simple models}

\label{section_polynomials}

We will use the results of the previous section to prove Theorems \ref%
{theorem_lin_statistics_Wigner} and \ref{theorem_lin_statistics_Wishart}
below, which are versions of Theorems \ref{theorem_lin_statistics_Wigner_general} and \ref{theorem_lin_statistics_Wishart_general} for models with a special distribution of entries.

\subsection{Wigner matrices}

 Assume in this section that $X$ is an infinite random Hermitian matrix with zero  
diagonal entries. We consider two possibilities for off-diagonal entries.
Either the entries take
values $\pm 1$ with probability $1/2,$ or they are  uniformly distributed on the unit circle. (The second model is not needed for the proof of Theorem \ref{theorem_lin_statistics_Wigner_general}. However, the results for this model can be obtained without any additional effort and they give a hint what happens in the situation with complex-valued entries.) 

Matrices  $A_N$ and $B_N$ are principal square submatrices of $X$ of the size $a^{\left(
N\right)}$ and $b^{\left(N\right) }$, respectively. The normalized matrices $\widetilde{A}_{N}$ are defined
as follows: 
\begin{equation*}
\widetilde{A}_{N}:=\frac{1}{2\sqrt{\left( a^{\left( N\right) }-2\right) }}%
A_{N},\text{ and }\widetilde{B}_{N}:=\frac{1}{2\sqrt{\left( b^{\left(
N\right) }-2\right) }}B_{N}.
\end{equation*}%
The choice of $a^{\left( N\right) }-2$ instead of $a^{\left( N\right) }$ and 
$b^{\left( N\right) }-2$ instead of $b^{\left( N\right) }$ is clearly not
essential for first-order asymtotics. However it makes some formulas
shorter. Recall that if $f:\mathbb{R}\rightarrow \mathbb{R}$ is a test
function, then its linear statistic for matrix $A_{N}$ is 
\begin{equation*}
\mathcal{N}\left( f,A_{N}\right) :=\mathrm{Tr}\left[ f\left( \widetilde{A}%
_{N}\right) \right] =\sum f\left( \lambda _{i}\right) ,
\end{equation*}%
where $\lambda _{i}$ are the eigenvalues of the matrix $\widetilde{A}_{N}.$
The quantity $\mathcal{N}\left( f,B_{N}\right) $ is defined similarly. The
centered statistics $\mathcal{N}^{o}\left( T_{k},A_{N}\right) $ and $%
\mathcal{N}^{o}\left( T_{l},B_{N}\right) $ are obtained by subtracting the
corresponding expectation values.

\begin{theorem}
\label{theorem_lin_statistics_Wigner} Let $T_{k}\left( x\right) $ denote the
Chebyshev polynomials of the first kind. As $N\rightarrow \infty ,$ the
distribution of the centered linear statistics $\mathcal{N}^{o}\left(
T_{k},A_{N}\right) $ and $\mathcal{N}^{o}\left( T_{l},B_{N}\right) $
converges to a two-variate Gaussian distribution with the covariance equal
to 
\begin{equation}
\left\{ 
\begin{array}{cc}
\delta _{kl}\frac{k}{2\beta }\left( \frac{\Delta }{\sqrt{ab}}\right) ^{k}, & 
\text{if }\min \left\{ k,l\right\} \geq 3, \\ 
0, & \text{if }\min \left\{ k,l\right\} \leq 2,%
\end{array}%
\right.  \label{formula_covariance0}
\end{equation}%
with $\beta =1$ for the model with $\pm 1$ entries and $\beta =2$ for the
model with entries on the unit circle.
\end{theorem}

\textbf{Proof of Theorem \ref{theorem_lin_statistics_Wigner}:} We will only
compute the limit covariance of $\mathcal{N}^{o}\left(
T_{k},A_{N}\right) $ and $\mathcal{N}^{o}\left( T_{l},B_{N}\right) $. The proof for higher moments follows from a similar combinatorial analysis. This analysis is
sketched below in the proof of the corresponding theorem for sample
covariance matrices.

Note that $A_N$ and $B_N$ are generalized adjacency matrices for complete graphs $G_{A_N}$ and $G_{B_N}$. 
These graphs are $(a_N-1)$- and $(b_N-1)$-regular, respectively. Therefore, by Theorem \ref{theorem_NBT_T_Wigner} the covariance of $\mathcal{N}^{o}\left(
T_{k},A_{N}\right) $ and $\mathcal{N}^{o}\left( T_{l},B_{N}\right) $ equals to: 
\begin{equation}
\frac{1}{4(a_N-2) ^{k/2}(b_N-2) ^{l/2}}
\sum_{\gamma,\gamma^{\prime}} 
\left[\E A(\gamma) B(\gamma^{\prime}) - \E A(\gamma) \E B (\gamma^{\prime}) \right] ,
\label{covariance_Wigner}
\end{equation}%
where the sums are over all pairs of non-backtracking tailless
(``NBT'') cyclic paths $\gamma$ and $\gamma^{\prime}$ of length $k$ and $l$, respectively, and $B (\gamma^{\prime})$ is defined analogously to  $A(\gamma)$.
Consider the sum 
\begin{equation}
\sum_{\gamma,\gamma^{\prime}} 
\left[\E A(\gamma) B(\gamma^{\prime}) - \E A(\gamma) \E B (\gamma^{\prime}) \right],
\label{second_moment}
\end{equation}%
and assign a graph and a type to each term in this sum. The graph is
formed by the edges of the NBT closed paths $\gamma$ and $\gamma^{\prime}$.   The \emph{type%
} of a pair of paths ($\gamma$, $\gamma^{\prime}$) is the graph together with a pair of paths on this graph, which are induced by $\gamma$ and $\gamma^{\prime}$. We
understand that the original labels of the vertices are removed in the sense that two pairs of paths $(\gamma,\gamma^{\prime})$ belong to the same type if they can be obtained each from the other by re-labeling of vertices. 

Matrices $A$ and $B$ are submatrices of an Hermitian matrix $X$ and by assumption the entries of $X$ are either $\pm 1$ with equal probability (case $\beta=1$) or uniformly distributed on the unit circle ($\beta=2$). It follows that
\begin{equation*}
\E \left[(A_{ij})^{k_1} (A_{ji})^{k_2} (B_{ij})^{k_3} (B_{ji})^{k_4}\right]=
\left \{
\begin{array}{ll}
1, & \text{if } \beta=1  \text{ and } k_1+k_2+k_3+k_4 \equiv 0 (\mathrm{mod}\ 2),\\
1, & \text{if } \beta=2 \text{ and } k_1+k_3=k_2+k_4 ,\\
0, & \text{otherwise.}
\end{array}
\right .
\end{equation*}
  In addition, the entries corresponding to different edges are independent.  Hence, the only types that bring a non-zero
contribution to the sum (\ref{second_moment}) are those in which paths traverse every edge an even
number of times.

Also, the contribution of every term with a disconnected graph is zero, since in this case the paths $\gamma$ and $\gamma^{\prime}$ are disjoint, and the contribution of $\E A(\gamma) B(\gamma^{\prime})$
is cancelled by the contribution of  $\E A(\gamma) \E B (\gamma^{\prime})$.

Now, consider the sum of all terms that have their graphs equal to the
cycle on $k$ vertices. This only happens if $l=k>2$ (there are no cycles of length $1$ or $2$ in the underlying graph). 

If $k=l>2$, then $\E A(\gamma)\E B(\gamma^{\prime})=0$ and $\E[A(\gamma)B(\gamma^{\prime})]=1$. (In the case $\beta=2$, this expectation is $1$ only if $\gamma$ and $\gamma^{\prime}$ have opposite orientations.) 

Note that each cyclic graph corresponds to $2k$ different types of paths $(\gamma,\gamma^{\prime})$ since every NBT
path can start from any of $k$ possible vertices in $V$ and the paths can
have either the same or opposite orientations. (Recall that the types are different only if they cannot be obtained by a re-labeling of the graph. Hence, the type is determined by the offset of the starting point and the orientation of one path relative to the other.) In the case $\beta =1$ both
orientations contribute and in the case $\beta =2,$ only opposite orientations
contribute. So the number of contributing types is $\frac{2}{\beta}k.$

In order to estimate the number of terms within each type, we note that since both paths are on the cycle, the
vertex labels are common to both $A$ and $B$ matrices. Since the number of rows and columns common to matrices $A$ and $B$ is by definition $\Delta^{(N)}$, the number of choices of $k$ (non-equal) labels is 
\begin{equation*}
(\Delta ^{(N)})(\Delta ^{(N)}-1)\ldots(\Delta ^{(N)}-k+1)=\left(\Delta ^{(N)}\right)^k+O\left(k^2\left(\Delta ^{(N)}\right)^{k-1}\right).
\end{equation*}    Hence, the total number of contributing terms with the cyclic graph can be
estimated as 
\begin{equation*}
\frac{2}{\beta}k\left(\Delta ^{(N)}\right)^k+O\left(k^3\left(\Delta ^{(N)}\right)^{k-1}\right).
\end{equation*}%
Each of these terms brings a contribution of $1$ to the sum (\ref{second_moment}).
 After the normalization given in (\ref{covariance_Wigner}), we find that
the contribution of these terms to the covariance is asymptotically $\frac{k}{2\beta}
\left(\frac{\Delta}{\sqrt{ab}}\right) ^k.$

The next step is to show that the contribution of all other terms is
negligible for all other types if $N$ is large. The crucial observation here is that every vertex in the graph of a pair 
($\gamma$,$\gamma^{\prime}$) has the degree greater or equal to two, since the paths $\gamma$ and $\gamma^{\prime}$ are NBT.

Consider the case $l=k$. Fix a type of ($\gamma$,$\gamma^{\prime}$) and assume that the graph is not a cycle of length $k$. First, suppose that there exists an edge of the graph which is traversed by the NBT paths more than twice. Then the
total number of edges in the graph is $<k,$ hence the number of
vertices is also less than $k$ (using the fact that each vertex must have
the degree of at least $2$). The number of types with this graph is bounded by a function of $k$ that counts the number of pairs of paths on the graph. The number of labellings of a graph with less than $k$ vertices is bounded by $(a^{(N)}+b^{(N)})^{k-1}$. The number of non-isomorphic graphs with less than $k$ vertices is bounded by another function of $k$. Hence the total number of pairs ($\gamma$,$\gamma^{\prime}$) with non-cyclic graphs is bounded by $ c_k(a^{(N)}+b^{(N)})^{k-1} $, where $c_k$ is a constant that depends on $k$ only. Finally, the contribution of each term is bounded by 1 by our assumption on the entries of the matrix $X$. Since we normalize by dividing by a multiple of $(t_N)^k$, these considerations imply that the sum of all considered terms gives a negligible contribution for
large $N$ provided that $k$ is fixed.

Next, suppose that every edge in a type's graph is traversed exactly twice by the NBT paths,
hence the number of edges is $k.$ Since the graph is not a
cycle, one of the vertices must have the degree of at least $3$
(since the graph is connected and all vertices have the degree of
at least $2$). Since the sum of vertex degrees is twice the number of edges it follows that the number of vertices is $<k$. By the same argument as in the previous paragraph the contribution of the sum of these terms is negligible.

If $k\neq l,$ then either the graph is not a cycle or some of the edges is traversed more than twice, and we find by a similar argument that contributions of all types are negligible. \hfill $\square $

\subsection{Sample covariance matrices}

In this section we assume that $A_{N}$ and $B_{N}$ are submatrices of an infinite random matrix $X$ whose entries are either $\pm 1$ with equal probability or are uniformly distributed on the unit circle. For the first case we set $\beta =1$ and for the second, $\beta =2$.

Let us define the normalized sample covariance matrices in the following
form: 
\begin{equation*}
W_{A_{N}}:=\frac{1}{2\sqrt{\left( a_{1}^{\left( N\right) }-1\right) \left(
a_{2}^{(N)}-1\right) }}\left[ A_{N}A_{N}^{\ast }-\left( a_{1}^{\left(
N\right) }+a_{2}^{\left( N\right) }-2\right) I_{a_{1}^{\left( N\right) }}%
\right] ,
\end{equation*}%
where $I_{a_{1}^{\left( N\right) }}$ is the $a_{1}^{\left( N\right) }$-by-$%
a_{1}^{\left( N\right) }$ identity matrix, and define $W_{B_{N}}$ similarly.

The normalization is chosen in such a way that the empirical distribution of
eigenvalues of $W_{A_{N}}$ and $W_{B_{N}}$ converges to a distribution
supported on the interval $\left[ -1,1\right] .$  Define 
\begin{equation*}
\mathcal{N}\left( f,W_{A_{N}}\right) :=\mathrm{Tr}\left[ f\left(
W_{A_{N}}\right) \right] =\sum f\left( \lambda _{i}\right) ,
\end{equation*}%
where $\lambda _{i}$ are the eigenvalues of the matrix $W_{A_{N}},$ and let 
\begin{equation}
\mathcal{N}^{o}\left( f,W_{A_{N}}\right) :=\mathcal{N}\left(
f,W_{A_{N}}\right) -\mathbb{E}\mathcal{N}\left( f,W_{A_{N}}\right) .
\label{cent_stat_Wishart_2}
\end{equation}%
The linear statistics for the matrix $W_{B_{N}}$ are defined similarly. 
Recall that there is a parameter $t_{N}$ that approaches infinity as $%
N\rightarrow \infty ,$ \ and that quantities $a_{1}^{(N)}/t_{N},$ $%
a_{2}^{\left( N\right) }/t_{N},$ $b_{1}^{\left( N\right) }/t_{N},$ $%
b_{2}^{\left( N\right) }/t_{N},$ $\Delta _{1}^{\left( N\right) }/t_{N},$ and 
$\Delta _{2}^{\left( N\right) }/t_{N}$ approach some positive limits $a_{1},$
$a_{2},$ $b_{1},$ $b_{2},$ $\Delta _{1},$ and $\Delta _{2},$ respectively.

\begin{theorem}
\label{theorem_lin_statistics_Wishart} Let $T_{k}\left( x\right) $ denote
the Chebyshev polynomials of the first kind. As $N\rightarrow \infty ,$ the
distribution of the centered linear statistics $\mathcal{N}^{o}\left(
T_{k},W_{A_{N}}\right) $ and $\mathcal{N}^{o}\left( T_{l},W_{B_{N}}\right) $
converges to a two-variate Gaussian distribution with the covariance equal
to 
\begin{equation*}
\delta _{kl}\frac{k}{2\beta }\gamma ^{k},
\end{equation*}%
if $\min \left\{ k,l\right\} >1$ and $0$ otherwise. Here, $\gamma :=\frac{\Delta _{1}\Delta _{2}}{\sqrt{a_{1}a_{2}b_{1}b_{2}}}$, $\beta =1$ for the
model with $\pm 1$ entries and $\beta =2$ for the model with entries on the
unit circle.
\end{theorem}

\textsc{Proof of Theorem \ref{theorem_lin_statistics_Wishart}:} 
Let
\begin{eqnarray} 
X_{k}\left( A_{N}\right)& =&\left(
\frac{\left( a_{1}^{\left( N\right)}-1\right) \left( a_{2}^{\left( N\right) }-1\right)}
{t_N^2} \right) ^{k/2}\mathcal{N}%
^{o}\left( T_{k},W_{A_{N}}\right) \notag \\
&=&\left( t_N\right) ^{-k}\left[ 
\mathcal{N}\left(\widetilde{T}_{k}[a_{1}^{\left( N\right) }-1,a_{2}^{\left( N\right)
}-1],A_{N}A_{N}^{\ast }\right) \right . \notag \\ 
&&-\left. \mathbb{E}\mathcal{N}\left( \widetilde{T} _{k}[a_{1}^{\left( N\right) }-1,a_{2}^{\left( N\right) }-1],A_{N}A_{N}^{\ast
}\right) \right]
\end{eqnarray}
where $\widetilde{T} _{k}$ are as in (\ref{definition_T_modified}): 
\begin{equation*}
\widetilde{T}_{k}[c,d]\left( x\right) :=\left( cd\right) ^{k/2}T_{k}\left( 
\frac{x-\left( c+d\right) }{2\sqrt{cd}}\right) ,
\end{equation*}
Let $X_{k}\left( B_{N}\right)$ be defined similarly.

Then it is enough to check that for $N\rightarrow \infty ,$ the
distribution of the random variables $X_{k}\left( A_{N}\right) $
and $X_{l}\left( B_{N}\right) $ converges to a two-variate Gaussian
distribution with the covariance matrix equal to 
\begin{equation}
\delta _{kl}\frac{k}{\beta }\left( 
\begin{array}{cc}
\left( a_{1}a_{2}\right) ^{k} & \left( \Delta _{1}\Delta _{2}\right) ^{k} \\ 
\left( \Delta _{1}\Delta _{2}\right) ^{k} & \left( b_{1}b_{2}\right) ^{k}%
\end{array}%
\right) ,
\label{covariance_Wishart}
\end{equation}%
if $\min \left\{ k,l\right\} >1$ and $0$ otherwise.

In order to prove this, note
that the matrix $A_{N}$ is a generalized adjacency matrix for a complete
bipartite graph $G$. The vertex sets of $G$ have $a_{1}^{\left( N\right) }$
and $a_{2}^{\left( N\right) }$ vertices, respectively. In particular, the
graph is $\left( a_{1}^{\left( N\right) },a_{2}^{\left( N\right) }\right) $%
-regular. By Theorem \ref{theorem_NBT_T_Wishart}, we see that 
\begin{equation*}
2X_k(A_N) =\frac{1}{%
\left( t_{N}\right) ^{k}}\sum A_2(\gamma) -r_A^{(N)},
\end{equation*}%
where the sum is over all NBT paths $\gamma$ that start with a vertex in $V$, and $r_A^{(N)}$ depends only on
$a_1^{(N)}$ and $a_2^{(N)}$ and, therefore, is non-random.
A similar formula holds for $X_k(B_N)$.

Since $r_{A}^{\left( N\right) }$ and $r_{B}^{\left( N\right) }$ are not
random, we only need to understand the joint distribution of  $\left( t_{N}\right) ^{-k}\sum A_2(\gamma) $ and $\left( t_{N}\right) ^{-l}\sum B_2(\gamma^{\prime})$, where $B_2(\gamma^{\prime})$ is defined analogously to $A_2(\gamma)$. In
particular, we need to show that in the limit of large $N$ this distribution
is Gaussian and to compute its covariance. The pattern of the argument is
well-known (see, for example, Section 2.1 in Anderson, Guionnet, Zeitouni 
\cite{anderson_guionnet_zeitouni10}) and for this reason we will be concise.

Consider the case $k=l.$ (The case of $k\neq l$ is similar and will be
omitted.)

\textbf{1. Covariance}
We are interested in estimating the following object: 
\begin{eqnarray}
\lim_{N\to \infty }\mathrm{Cov}\left( X_{k}\left( A_{N}\right)
,X_{k}\left( B_{N}\right) \right) \notag \\
=\lim_{N\to \infty }\frac{1}{4\left( t_{N}\right) ^{2k}} \sum_{\gamma,\gamma^{\prime}} \left[ \E A_2 (\gamma) \E B_2(\gamma^{\prime}) 
-\E A_2 (\gamma) \E B_2(\gamma^{\prime}) \right] , \notag
\end{eqnarray}%
where the sum is over all pairs of non-backtracking tailless
(``NBT'' ) cyclic paths of length $2k$ each.
Note that for $k=1,$ the number of such paths is zero. In the following we
assume $k>1$. 

Consider the normalized sum 
\begin{equation*}
 \frac{1}{\left( t_{N}\right) ^{2k}}\sum_{\gamma,\gamma^{\prime}} \left[ \E A_2 (\gamma) \E B_2(\gamma^{\prime}) 
-\E A_2 (\gamma) \E B_2(\gamma^{\prime}) \right],
\end{equation*}%
and assign a graph and a type to each term in this sum. This is done
as for Wigner matrices with a small modification. Namely, since the original
graph is bipartite, the term graphs are also bipartite and we will
keep the information about the partition.

Let $v \in V$ and $w \in W$. Consider the term where path $\gamma$ goes $k_1$ times from $v$ to $w$, and $k_2$ times from $w$ to $v$. In addition, let $\gamma^{\prime}$ go $k_3$ times from $v$ to $w$, and $k_4$ times from $w$ to $v$. 
Matrices $A$ and $B$ are submatrices of matrix $X$ and the entries of $X$ are either $\pm 1$ with equal probability (case $\beta=1$) or uniformly distributed on the unit circle ($\beta=2$).  Therefore, 
\begin{equation*}
\E \left[A_{vw}^{k_1}\overline{A_{vw}}^{k_2} B_{vw}^{k_3}\overline{B_{vw}}^{k_4}\right]=
\left \{
\begin{array}{ll}
1, & \text{if } \beta=1  \text{ and } k_1+k_2+k_3+k_4 \equiv 0 (\mathrm{mod}\ 2),\\
1, & \text{if } \beta=2 \text{ and } k_1+k_3=k_2+k_4 ,\\
0, & \text{otherwise.}
\end{array}
\right .
\end{equation*}
In addition, the entries corresponding to different edges are independent.  Hence, the only types that bring a non-zero
contribution to the sum are those in which paths traverse every edge an even
number of times. 

Similar to the case of Wigner matrices, the contribution of every term with a
disconnected graph is zero.

Now, consider the sum of all terms with the graph equal to
the cycle on $2k$ vertices. Since both paths are on this cycle, the
vertices must have the labels that are common to both $%
A$ and $B$ matrices. Repeating the argument from the proof of Theorem 3.1, we find that the contribution of the terms with 
the cyclic graph has the limit
 ($2/\beta ) k\left( \Delta _1 \Delta _2\right)^k. $ (The first factor is from the choice of path orientations on the cycle, and the second factor is from the choice of starting points. Note that although the cycle has the length $2k$, there are only $k$ possible starting points since the paths must start with a vertex in a particular partition.)  

The next step is to show that the contribution of all other types is
negligible if $N$ is large. This is done in the same way as in the proof of Theorem 3.1. The crucial observation is that the paths $\gamma$ and $\gamma^{\prime}$ are NBT and therefore every
vertex in a graph of the pair ($\gamma$,$\gamma^{\prime}$) must have the degree greater or equal to
two.

 We conclude that 
\begin{equation*}
\lim_{N\rightarrow \infty }\mathrm{Cov}\left( X_{k}\left( A_{N}\right)
,X_{k}\left( B_{N}\right) \right) =\frac{k}{2\beta }\left( \Delta _{1}\Delta _{2}\right) ^{k}.
\end{equation*}%
Similarly, 
\begin{eqnarray*}
\lim_{N\rightarrow \infty }\mathrm{Cov}\left( X_{k}\left( A_{N}\right)
,X_{k}\left( A_{N}\right) \right) &=&\frac{k}{2\beta }\left(
a_{1}a_{2}\right) ^{k},\text{ and } \\
\lim_{N\rightarrow \infty }\mathrm{Cov}\left( X_{k}\left( B_{N}\right)
,X_{k}\left( B_{N}\right) \right) &=&\frac{k}{2\beta }\left(
b_{1}b_{2}\right) ^{k}.
\end{eqnarray*}

\textbf{2. Higher moments}

The argument for higher moments is similar. Consider an $m$-th joint moment, 
\begin{equation*}
\mathbb{E}\left( \left[ \left(2 t_{N}\right) ^{-k}\sum A_2(\gamma) \right] ^{m_{a}}\left[
\left( 2 t_{N}\right) ^{-l}\sum B_2(\gamma^{\prime}) \right] ^{m_{b}}\right) ,
\end{equation*}%
where $m_{a}+m_{b}=m$ and $\gamma$ and $\gamma^{\prime}$ are NBT cyclic paths of length $2k$ and $2l$, respectively.  When we expand this expression, we obtain, 
\begin{equation*}
\left(2 t_{N}\right) ^{-km_a-lm_b} \sum \E A_2(\gamma_1) \ldots A_2(\gamma_{m_a})B_2(\gamma_{m_a+1}) \ldots B_2(\gamma_m) ,
\label{expansion_high_moments}
\end{equation*}% 
The type of a term in this expansion is
given by a graph and $m$ NBT cyclic paths. Non-negligible
contributions arise only if $k=l$ and come from the types whose graph is the union of
disjoint cycles of length $2k$. Every edge in these cycles must be traversed by the NBT
paths exactly twice, and therefore exactly two paths traverse a cycle. Therefore, every type with a non-negligible contribution corresponds to a matching on the set of $m$ NBT paths: The paths are matched if they traverse the same cycle. 

Let a particular matching on paths be fixed. The number of the $m$-tuples of paths that correspond to this matching converges asymptotically to the product $c_1 c_2 \ldots c_{m/2}\left( t_{N}\right) ^{km}$, where $c_i=\frac{2k}{\beta }\left(
a_1 a_2\right) ^{k},$ if the $i$-th pair in the matching pairs $\gamma_s$ with $\gamma_t$, and $\max\{s,t\} \leq m_a$ (that is, the paths in the $i$-th pair are both from the factors $A_2(\gamma)$ in \ref{expansion_high_moments}), 
$c_i=\frac{2k}{\beta }\left(b_1 b_2\right) ^{k},$ if $\min\{s,t\} > m_a$, and $c_i=\frac{2k}{\beta }\left(\Delta_1 \Delta_2\right) ^{k},$ otherwise. 

 Hence, in the limit the
sum of normalized contributions over all matchings coincides exactly with the Wick formula for the higher moments of the two-variate Gaussian distribution with the covariance matrix (\ref{covariance_Wishart}). (See Zee \cite{zee03}). \hfill $\square $

\section{Matrices with more general distribution of matrix entries}

\label{section_general_matrices}

\subsection{Wigner}

\textbf{Proof of Theorem \ref{theorem_lin_statistics_Wigner_general}}: 
The idea of the proof is to use Wigner matrices from Section \ref{section_polynomials}, for which the covariances
of the Chebyshev polynomials have been already computed, and then calculate how a change in the moments of matrix entries affect these covariances.  
  
The key formula is the following expression, valid for matrices that satisfy assumptions of
 Theorem \ref{theorem_lin_statistics_Wigner_general}, 
\begin{eqnarray}
&&\lim_{N\to \infty} \Cov \left(\Tr\left[ \left( A_{N}/\sqrt{t_{N}}\right) ^{k}\right],
 \Tr\left[ \left( B_{N}/\sqrt{t_{N}}\right) ^{l}\right] \right) \notag = \\
&&d_{2}klC_{(k-1)/2}C_{\left( l-1\right) /2}\Delta a^{\frac{k-1}{2}}b^{\frac{%
l-1}{2}}
+\left( m_{4}-1\right) \frac{kl}{2} C_{k/2}C_{l/2}\Delta ^{2}a^{\frac{k}{2}%
-1}b^{\frac{l}{2}-1} \nonumber \\
&&+\sum_{r=3}^{\infty }\frac{2kl}{r}\left( \sum\limits_{\substack{ s_{i}\geq
0  \\ 2\sum s_{i}=k-r}}\prod C_{s_{i}}\right) \left( \sum\limits_{\substack{ %
t_{i}\geq 0  \\ 2\sum t_{i}=l-r}}\prod C_{t_{i}}\right) \Delta ^{r}a^{\frac{%
k-r}{2}}b^{\frac{l-r}{2}}. \label{borodins_formula}
\end{eqnarray}%
Here $C_{k}:=\binom{2k}{k}/\left( k+1\right) $ are the Catalan numbers, which count the number of planar rooted trees with $k$ edges. If $k$ is not integer, then we set $C_k=0$.

This formula comes from counting the contribution of various types in the expansion similar to (\ref{covariance_Wigner}). It is significantly more complicated because the paths can now be backtracking and can have loops associated to the diagonal terms. The first term comes from the contribution of the paths ($\gamma$, $\gamma'$) that traverse two trees with $(k-1)/2$ and $(l-1)/2$ edges, respectively. These trees are disjoint except that they hang from a common vertex and the paths have a loop at this vertex. The factor $kl$ comes from a choice of starting points on the trees and the factor $\Delta a^{\frac{k-1}{2}}b^{\frac{l-1}{2}}$ comes from the label counting. (In particular, $\Delta$ comes from the number of labels for the common vertex.)  The factor $d_2$ also comes from the common vertex. 

Similarly, the second term comes from two trees with $k/2$ and $l/2$ edges, respectively, that are glued along one edge. There are $k/2 \times l/2$ choices for this edge and 2 possible orientations for the gluing. The factor $m_4-1$ comes from the glued edge. 

The third term comes from two graphs each of which is a cycle of length $r$ with attached trees. The graphs are glued along the cycle. The third term equals to the count of such graphs, multiplied by $kl$, which is the number of choices of starting points, and by $2$, which is the number of possible orientations for the gluing.  

One can find a sketch of the proof for this formula (for the case $d_2=2$, $m_4=3$) in Borodin's paper \cite{borodin10} in the discussion after formula (4). We omit the details.

We can understand the limit covariance of linear statistics as a bi-linear functional on the space of all polynomials:
\begin{equation}
\alpha(P,Q):=\lim_{N\to \infty} \Cov \left(\Tr\left[ P\left( \frac{A_N}{\sqrt{t_{N}}}\right)\right],
 \Tr\left[ Q\left(\frac{B_N}{\sqrt{t_{N}}}\right)\right] \right),
\end{equation}
where $P$ and $Q$ are two polynomials.

Then, (\ref{borodins_formula}) gives a representation of this bi-linear functional as a sum of three bi-linear functionals (corresponding to the three terms in the formula (\ref{borodins_formula})): 
\begin{equation}
\alpha(P,Q)=\alpha_1(P,Q)+\alpha_2(P,Q)+\alpha_3(P,Q).
\label{alpha_sum}
\end{equation}
We will evaluate these functionals at the Chebyshev's polynomials $\widetilde{T}_{k,a}\left( x\right) :=T_{k}\left( x/2%
\sqrt{a}\right) $.

Note that by (\ref{borodins_formula}), $\alpha_1(x^k,x^l)$ is not zero only if both $k$ and $l$ are odd, in which case,
\begin{eqnarray}
\alpha_1(x^k,x^l)&=&d_2 k l C_{(k-1)/2} C_{( l-1)/2} \Delta a^{\frac{k-1}{2}}b^{\frac{l-1}{2}} \notag\\
&=&d_2 \binom{k}{\frac{k-1}{2}}\binom{l}{\frac{l-1}{2}}\Delta a^{\frac{k-1}{2}}b^{\frac{l-1}{2}} \notag\\
&=&\frac{d_2 \Delta}{( 2\pi i)^2}\int\limits_{|z|=c_1}\left( z+\frac{a}{z}\right) ^k\frac{dz}{z^2 } \notag \\
&&\times
\int\limits_{|w|=c_2}\left( w+\frac{b}{w}\right) ^l\frac{dw}{w^2}. \label{contour_Wigner}
\end{eqnarray}

By the standard property of the Chebyshev's polynomials,  
$T_k( \frac{1}{2} (x + x^{-1})) = \frac{1}{2} (x^k + x^{-k})$ for $x \ne 0$. Hence,
\begin{equation}
\widetilde{T}_{k,a}\left( z+\frac{a}{z}\right) 
=T_k\left(\frac{1}{2} \left(\frac{z}{\sqrt{a}}+\frac{\sqrt{a}}{z}\right) \right)
=\frac{1}{2}\left[ \left( 
\frac{z}{\sqrt{a}}\right) ^{k}+\left( \frac{\sqrt{a}}{z}\right) ^{k}\right] .
\label{Chebyshevs_property}
\end{equation}%

By (\ref{contour_Wigner}), (\ref{Chebyshevs_property}), and the bi-linearity of $\alpha_1$, we have:

\begin{eqnarray*}
\alpha_1(\widetilde{T}_{k,a},\widetilde{T}_{l,b})
&=&\frac{d_2 \Delta}{\left( 2\pi i\right) ^{2}}
\int_{\vert z \vert =c_1}\widetilde{T}_{k,a}\left( z +\frac{a}{z}\right)\frac{dz}{z^2} \\
&&\times \int_{\left\vert w\right\vert =c_{2}}
\widetilde{T}_{l,b} \left( w +\frac{b}{w}\right) \frac{dw}{w^2}\\
&=&\frac{d_2 \Delta /4}{\left( 2\pi i\right) ^{2}}
\int_{|z| =c_{1}}\left[ \left( \frac{z}{\sqrt{a}}\right) ^{k}+\left( 
\frac{\sqrt{a}}{z}\right) ^{k}\right] \frac{dz}{z^2} \\
&&\times \int_{\left\vert w\right\vert =c_{2}}\left[ \left( \frac{w}{\sqrt{b}%
}\right) ^{l}+\left( \frac{\sqrt{b}}{w}\right) ^{l}\right] \frac{dw}{w^2}.
\end{eqnarray*}%

This is different from zero if and only if $k=l=1,$ in which case we find
\begin{equation}
\alpha_1(\widetilde{T}_{1,a},\widetilde{T}_{1,b})=\frac{d_{2}}{4}\frac{\Delta }{\sqrt{ab}}.
\label{first_alpha} 
\end{equation}%

Similarly, $\alpha_2(x^k,x^l)$ is non-zero only if both $k$ and $l$ are even, and then,
\begin{eqnarray}
\alpha_2(x^k,x^l)&=&\left( m_{4}-1\right) \frac{kl}{2} C_{k/2}C_{l/2}\Delta ^{2}a^{\frac{k}{2}%
-1}b^{\frac{l}{2}-1} \notag \\
&=&2( m_4-1)\binom{k}{\frac{k}{2}-1}\binom{l}{\frac{l}{2}-1}\Delta ^{2}a^{\frac{k}{2}-1}b^{\frac{l}{2}-1} \notag \\
&=&\frac{2(m_4-1) \Delta^2}{( 2\pi i)^2}\int\limits_{|z|=c_1}\left( z+\frac{a}{z}\right) ^k\frac{dz}{z^3} \notag \\
&&\times
\int\limits_{|w|=c_2}\left( w+\frac{b}{w}\right) ^l\frac{dw}{w^3}. \label{contour_Wigner2}
\end{eqnarray}

Therefore,

\begin{eqnarray*}
\alpha_2(\widetilde{T}_{k,a},\widetilde{T}_{l,b})
&=&\frac{2(m_4-1) \Delta^2}{( 2\pi i)^2}
\int_{\vert z \vert =c_1}\widetilde{T}_{k,a}\left( z +\frac{a}{z}\right) \frac{dz}{z^3} \\
&&\times \int_{\left\vert w\right\vert =c_{2}}
\widetilde{T}_{l,b} \left( w +\frac{b}{w}\right) \frac{dw}{w^3}\\
&=&\frac{(m_4-1) \Delta^2 /2}{( 2\pi i)^2}
\int_{|z| =c_{1}}\left[ \left( \frac{z}{\sqrt{a}}\right) ^{k}+\left( 
\frac{\sqrt{a}}{z}\right) ^{k}\right] \frac{dz}{z^3} \\
&&\times \int_{\left\vert w\right\vert =c_{2}}\left[ \left( \frac{w}{\sqrt{b}%
}\right) ^{l}+\left( \frac{\sqrt{b}}{w}\right) ^{l}\right] \frac{dw}{w^3}.
\end{eqnarray*}%

This expression is not zero if and only if $k=l=2,$ in which
case it is 
\begin{equation}
\alpha_2(\widetilde{T}_{2,a},\widetilde{T}_{2,b})=\frac{m_{4}-1}{2}\frac{\Delta ^{2}}{ab}.
\label{second_alpha}
\end{equation}

Suppose now we take $A_N$ and $B_N$ in equation (\ref{borodins_formula}) to be
Wigner matrices in Section \ref{section_polynomials}. In this case $d_2 = m_4 - 1 = 0$. Theorem \ref{theorem_lin_statistics_Wigner} and
equation (\ref{borodins_formula}) then give

\begin{equation*}
\alpha_3(\widetilde{T}_{k,a},\widetilde{T}_{l,b})=
\left\{ 
\begin{array}{cc}
\delta _{kl}\frac{k}{2}\left( \frac{\Delta }{\sqrt{ab}}\right) ^{k}, & 
\text{if }\min \left\{ k,l\right\} \geq 3, \\ 
0, & \text{if }\min \left\{ k,l\right\} \leq 2.%
\end{array}%
\right.  
\end{equation*}
This formula, and formulas (\ref{alpha_sum}), (\ref{first_alpha}), and (\ref{second_alpha}) complete the proof of Theorem \ref{theorem_lin_statistics_Wigner_general}. \hfill $\square$

\subsection{Sample Covariance}

We proceed as in the previous section about Wigner matrices. 
The covariance of the linear statistics for a special model
 has been already computed in Theorem \ref{theorem_lin_statistics_Wishart}). We define the bi-linear form

\begin{equation}
\alpha(P,Q;\left .\{m_k\}\right|_{k=3}^{\infty}):=\lim_{N\to \infty} \Cov \left(\Tr\left[ P\left( W_{A_N}\right)\right],
 \Tr\left[ Q\left(W_{B_N}\right)\right] \right),
\label{covariance_formula}
\end{equation}
where $P$ and $Q$ are two polynomials, and $m_k$ are the moments of the matrix entries. Let 
$\left .\{m'_k\}\right|_{k=3}^{\infty}$ denote a specific sequence of moments, with $m'_k=0$ for odd $k$ and $m'_k=1$ for even $k$. Define
\begin{eqnarray}
\alpha_2(P,Q)&:=&\alpha(P,Q,\left .\{m'_k\}\right|_{k=3}^{\infty}),\\
\alpha_1(P,Q)&:=&\alpha(P,Q,\left .\{m_k\}\right|_{k=3}^{\infty})-\alpha_2(P,Q).
\label{def_alpha1}
\end{eqnarray}

Then, 
\begin{equation}
\alpha(P,Q,\left .\{m_k\}\right|_{k=3}^{\infty})=\alpha_1(P,Q)+\alpha_2(P,Q).
\label{alpha_sum2}
\end{equation}
Since $\alpha_2(P,Q)$ is known from Theorem \ref{theorem_lin_statistics_Wishart}, we only need to calculate $\alpha_1$. 

The quantities $\alpha_1(x^k,x^l)$ are calculated in the next two results. 

\begin{lemma}
\label{prop_contribution_higher}
\begin{equation*}
\alpha_1(x^k,x^l)=\left( m_{4}-1\right) \gamma
R_{k}\left( a_{1},a_{2}\right) R_{l}\left( b_{1},b_{2}\right) .
\end{equation*}
\end{lemma}

Here, $R_{k}\left( a_{1},a_{2}\right)$, $k\geq 1,$ is defined as follows:%
\begin{eqnarray*}
R_{k}\left( a_{1},a_{2}\right) &=&\left( \frac{a_{1}+a_{2}}{2\sqrt{a_{1}a_{2}%
}}\right) ^{k}\sqrt{\frac{a_{2}}{a_{1}}} \\
&&\times \sum\limits_{s=1}^{k}\left( \frac{-1}{a_{1}+a_{2}}\right) ^{s}%
\binom{k}{s}s\sum_{t=1}^{s}\stack{s}{t}a_{1}^{t}a_{2}^{s-t},
\end{eqnarray*}%
where  
\begin{equation*}
\stack{s}{t}:=\frac{1}{s}\binom{s-1}{t-1}\binom{s}{t-1}.
\end{equation*}%
These coefficients are called Narayana numbers and they frequently occur in
combinatorial problems. In particular, they count the number of rooted
planar trees with $n$ edges and $k$ leaves, and also the number of partitions of
the set $\left[ n\right] =\{1,\ldots ,n\}$ that have $k$ blocks.

This lemma follows from the combinatorial analysis of graph contributions. Its proof is postponed to the end of this section.
Surprisingly, $R_{k}\left( a_{1},a_{2}\right) $ can be computed quite explicitly and its value does not depend on $a_{1}$ or $a_{2}.$

\begin{lemma}
\label{lemma_value_Rk}%
\begin{equation*}
R_{k}\left( a_{1},a_{2}\right) =\left\{ 
\begin{array}{cc}
-\frac{1}{2^{k}}\binom{k}{\frac{k-1}{2}}, & \text{if }k\text{ is odd,} \\ 
0 & \text{if }k\text{ is even.}%
\end{array}%
\right.
\end{equation*}
\end{lemma}
The proof of this lemma is also postponed to the end of this section.

\textbf{Proof of Theorem \ref{theorem_lin_statistics_Wishart_general}}: By
using Lemmas \ref{prop_contribution_higher} and Lemma \ref{lemma_value_Rk}, we find that $\alpha_1(x^k,x^l)$
is non-zero if and only if $k$ and $l$ are odd, and then,
\begin{eqnarray}
\alpha_1(x^k,x^l)&=&
\frac{( m_4-1) \gamma}{2^{k+l}}\binom{k}{\frac{k-1}{2}} \binom{l}{\frac{l-1}{2}} \notag \\
&=&\frac{(m_4-1) \gamma}{( 2\pi i)^2}\int\limits_{|z|=c_1}\left(\frac{1}{2}\left( z+\frac{1}{z}\right)\right) ^k\frac{dz}{z^2} \notag \\
&&\times
\int\limits_{|w|=c_2}\left(\frac{1}{2}\left( w+\frac{1}{w}\right)\right) ^l\frac{dw}{w^2}.
\end{eqnarray}

Hence,

\begin{eqnarray*}
\alpha_1(T_k,T_l)
&=&\frac{(m_4-1) \gamma}{( 2\pi i)^2}
\int_{\vert z \vert =c_1}T_k\left(\frac{1}{2}\left( z +\frac{1}{z}\right)\right) \frac{dz}{z^2} \\
&&\times \int_{\left\vert w\right\vert =c_{2}}
T_l \left(\frac{1}{2}\left( w +\frac{1}{w}\right)\right) \frac{dw}{w^2}\\
&=&\frac{(m_4-1) \gamma}{( 2\pi i)^2}
\int_{|z| =c_{1}}\frac{1}{2}\left[z^k+z^{-k}\right] \frac{dz}{z^2} \\
&&\times \int_{\left\vert w\right\vert =c_{2}}\frac{1}{2}\left[w^l+w^{-l}\right] \frac{dw}{w^2},
\end{eqnarray*}%
which is non-zero if and only if $k=l=1,$ in which case 
\begin{equation*}
\alpha_1(T_1,T_1)=\frac{\left( m_{4}-1\right) \gamma }{4}.
\end{equation*}%

Then, equation (\ref{alpha_sum2}) and Theorem \ref{theorem_lin_statistics_Wishart} give the statement of the theorem.
\hfill $\square $

\textbf{Proof of Lemma \ref{prop_contribution_higher}: } 
We use the definition of $W_{A_N}$ and $W_{B_N}$ in equations (\ref{W_A_definition}) and (\ref{W_B_definition}), and note that by the binomial theorem, 
\begin{eqnarray}
\label{formula_binomial}
\Tr\left( \frac{A_{N}A_{N}^{\ast }}{2\sqrt{a_{1}a_{2}}N}-\frac{%
a_{1}+a_{2}}{2\sqrt{a_{1}a_{2}}}I_{N}\right) ^{k} &=&\sum_{s=1}^{k}\left(
-1\right) ^{s}\binom{k}{s}\Tr\left( \frac{A_{N}A_{N}^{\ast }}{2\sqrt{%
a_{1}a_{2}}N}\right) ^{s}\left( \frac{a_{1}+a_{2}}{2\sqrt{a_{1}a_{2}}}%
\right) ^{k-s} \notag \\
&&+N\left( \frac{a_{1}+a_{2}}{2\sqrt{a_{1}a_{2}}}\right) ^{k},
\end{eqnarray}%
and similarly for matrix $B_{N}$.
Next, we consider the following expression: 
\begin{equation}
\E\left( \Tr\left( A_{N}A_{N}^{\ast }\right) ^{k}\Tr\left( B_{N}B_{N}^{\ast }\right) ^{l}\right) 
- \E\Tr\left( A_{N}A_{N}^{\ast }\right) ^{k}
\E\Tr\left( B_{N}B_{N}^{\ast }\right) ^{l}.
\label{expression_covariance}
\end{equation}

By formula (\ref{formula_binomial}) the covariance in definition (\ref{covariance_formula}) is a linear combination of these expressions for appropriate $k$ and $l$.%

As usual, we expand expression (\ref{expression_covariance}) as a sum of the expected products of matrix entries. Every product can be coded by two paths on a graph. The following statements are easy to check. Every edge in a graph must be traversed at least twice. Only connected graphs contribute. The only graphs that contribute are therefore either trees, or graphs with
only one cycle. The types corresponding to the graphs with a cycle contribute only if each edge is traversed exactly twice. Hence the contribution of such graphs involves only the second moment of matrix entries. This contribution is the same for the special model, and therefore (by definition 
(\ref{def_alpha1}) and by linearity) these graphs contribute zero to $\alpha_1(x^k,x^l)$. 

The graphs that involve moments higher than the second are trees. Each of the two paths on a tree traverse a sub-tree and the entire tree is given by these (two) sub-trees glued along an edge. One of these sub-trees has $k$ edges and another one has $l$ edges. The corresponding paths have lengths $2k$ and $2l,$ respectively, and they
traverse every edge of the corresponding tree twice.

Recall that the graphs are partitioned. Some of its vertices correspond to rows of matrices $A_N$ and $B_N$ and some to columns. In order to count valid labellings corresponding to a particular type, suppose that the sub-tree with $k$ edges has $t_{1}$ row vertices and $k+1-t_{1}$ column vertices, and the sub-tree with $l$ edges has $t_{2}$ row
vertices and $l+1-t_{2}$ column vertices. (We set the matrix size parameter $t_{N}$ equal to $N$ here, and use letter $t$ for a different purpose.) The subtree with $k$ vertices corresponds to matrix $A_N$  and its row and column vertices can be labeled in $a_{1}^{(N)}$ and $a_{2}^{(N)}$ different ways, respectively. The only exceptions are a row and a column vertices that belong to glued edge. These vertices can be labeled in $\Delta_1^{( N) }$ and $\Delta_2^{( N) }$ different ways, respectively. Similarly for the subtree that corresponds to matrix $B_N$ . Thus, the contribution of this type  is 
\begin{eqnarray}
( m_{4}-1)\Delta_1^{( N) } ( a_{1}^{(N)})
^{t_{1}-1}( a_{2}^{( N) }) ^{k-t_{1}}
\Delta_2^{( N) } ( b_1^{( N)})
^{t_{2}-1}( b_2^{( N)}) ^{l-t_{2}} \notag \\
+O\left( N^{k+l-1}\right) . \notag
\end{eqnarray}%
(Here, we have the factor $m_4-1$ and not $m_4$ because we subtract the corresponding contribution for the special model with $m'_4=1$.)
 Hence, for 
\begin{equation}
\Cov\left(\Tr\left( \frac{A_{N}A_{N}^{\ast }}{2%
\sqrt{a_{1}a_{2}}N}\right) ^{k},\Tr\left( \frac{B_{N}B_{N}^{\ast }}{2%
\sqrt{b_{1}b_{2}}N}\right) ^{l}\right) ,  \label{cov_noncentered_monomials}
\end{equation}%
as $N\rightarrow \infty ,$ the contribution of this type converges to 
\begin{eqnarray}
&&\left( m_{4}-1\right) \frac{1}{2^{k+l}}\left( \frac{a_{2}}{a_{1}}\right)
^{k/2-t_{1}}\left( \frac{b_{2}}{b_{1}}\right) ^{l/2-t_{2}}\frac{\Delta _{1}}{%
a_{1}}\frac{\Delta _{2}}{b_{1}}  \notag \\
&=&\left( m_{4}-1\right) \gamma \frac{1}{2^{k+l}}\left( \frac{a_{2}}{a_{1}}%
\right) ^{k/2-t_{1}+1/2}\left( \frac{b_{2}}{b_{1}}\right) ^{l/2-t_{2}+1/2}.
\label{contribution}
\end{eqnarray}

In order to calculate the number of types with this contribution, we need the
following lemma.

\begin{lemma}
The number of the non-isomorphic bipartite planar rooted trees with $n$
edges that have $k$ vertices in one of the partitions equals the Narayana
number $\stack{n}{k}.$
\label{graph_equivalence}
\end{lemma}
Let us postpone the proof of this lemma. It implies that the number of types with contribution (\ref{contribution}%
) is given by 
\begin{equation*}
kl\stack{k}{t_{1}}\stack{l}{t_{2}}.
\end{equation*}%
(The prefactor $kl$ corresponds to the choice of the tree edges that are
glued together in the graph.) Therefore, the total contribution of
trees to (\ref{cov_noncentered_monomials}) converges to 
\begin{equation}
\frac{\left( m_{4}-1\right) \gamma kl}{2^{k+l}}\sum_{t_{1}=1}^{k}
\stack{k}{t_{1}}\left( \frac{a_{2}}{a_{1}}\right) ^{k/2-t_{1}+1/2}\sum_{s_{2}=1}^{l}
\stack{l}{t_{2}}\left( \frac{b_{2}}{b_{1}}\right) ^{l/2-t_{2}+1/2}.
\label{cov_noncentered_monomials_value}
\end{equation}%
 
By using (\ref{formula_binomial}) and (\ref{cov_noncentered_monomials_value}), we obtain the formula in
the claim of Lemma \ref{prop_contribution_higher}. \hfill $\square $

\begin{figure}[tbph]
\includegraphics[width=6cm]{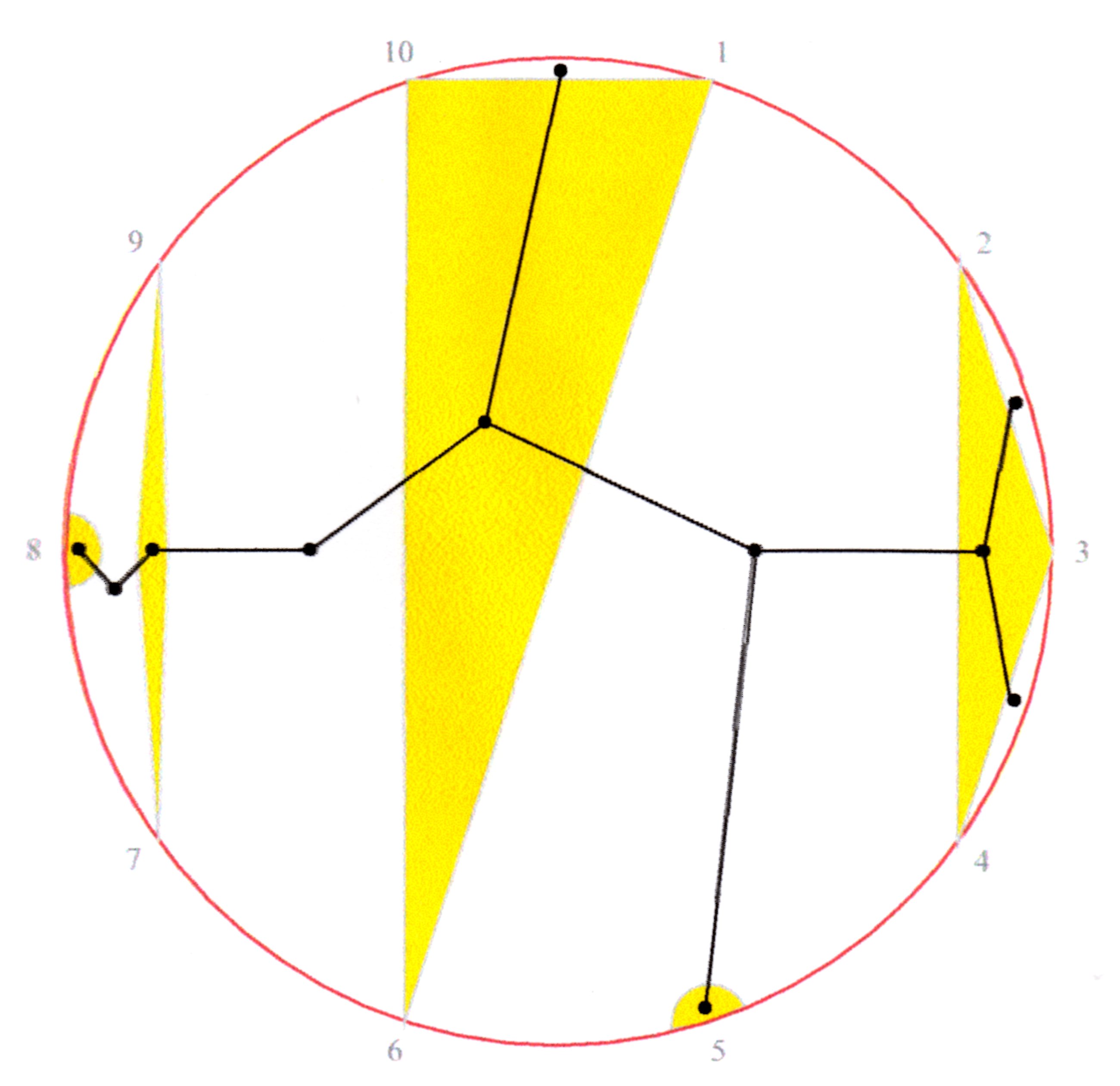} 
\caption{ An NC polygon diagram with superimposed bicolored plane tree.
Polygon sides $\leftrightarrow $ tree edges. The NC partition is $(1, 6,10)(2,3,4)(5)(7 9)(8)$.  (Graphics is courtesy Callan
and Smiley.)}
\label{fig:bijection}
\end{figure}

\textbf{Proof of Lemma \ref{graph_equivalence}:} 
The proof is based on a bijection between bipartite planar rooted trees and polygon systems corresponding to non-crossing partitions. The bijection appears in Callan and Smiley in \cite{callan_smiley05} who use it to derive several results about non-crossing partitions. However, the connection with Narayana numbers has not been noticed.

First, let us explain the correspondence between non-crossing partitions and polygonal systems. 
Recall that a non-crossing partition of $\left[ n\right] =\{1,\ldots ,n\}$
is one for which no $4$-tuple $a<b<c<d$ has $a$ and $c$ in one block and $b$
and $d$ in another. This implies that if the elements of $\left[ n\right] $
are realized as points on the circle, and neighboring elements within each
block are joined by line segments, then a non-crossing partition will appear
as a system of non-overlapping polygons. It is clear that the number of polygons equals to
the number of blocks in the partition. 

The bijection of these polygonal systems with bipartite planar trees works as follows. Let the polygons be colored yellow and the remaining regions of the disk colored white. Place a vertex in each region of the disk, both yellow and
white. Join vertices in adjacent regions by edges. Then allow each vertex to inherit the
color of the region it is in to get the desired bicolored plane tree. Figure \ref{fig:bijection} explains this bijection with an example. See \cite{callan_smiley05} for a proof that this is indeed a bijection.

Note that the number of
vertices in one of the partitions of the tree corresponds to the number of polygons in the polygonal system. It follows that the number of the bipartite planar rooted trees with $n$ edges and $k$ vertices
in one of the partitions equals the number of non-crossing partitions of 
$\left[ n\right] $ with $k$ blocks, which is known to be equal to the Narayana number $\stack{n}{k}$. \hfill $\square $

\textbf{Proof of Lemma \ref{lemma_value_Rk}:} Recall that the Narayana polynomials are defined as 
\begin{equation*}
N_{n}\left( x\right) =\sum_{k=1}^{n}\stack{n}{k}x^{k}.
\end{equation*}%
Hence, if we take $x=a_{1}/a_{2},$ then we have 
\begin{eqnarray*}
R_{k}\left( a_{1},a_{2}\right) &=&\frac{1}{2^{k}}\left( 1+x\right) ^{k}x^{-%
\frac{k+1}{2}} \\
&&\times \sum_{s=1}^{k}\left( \frac{-1}{1+x}\right) ^{s}s\binom{k}{s}%
N_{s}\left( x\right) .
\end{eqnarray*}

Next, we use the fact that $N_{s}\left( x\right) $ are related to a
particular case of the Jacobi polynomials. Namely, 
\begin{equation*}
N_{n}\left( x\right) =\frac{x}{n}\left( x-1\right) ^{n-1}P_{n-1}^{\left(
1,1\right) }\left( \frac{x+1}{x-1}\right) .
\end{equation*}%
(This fact was apparently first noted in Proposition 6 of \cite%
{kostov_mf_shapiro09}.) By substituting this identity into the previous
formula, we get 
\begin{eqnarray}
R_{k}\left( a_{1},a_{2}\right) &=&\frac{1}{2^{k}}\left( \frac{1+x}{\sqrt{x}}%
\right) ^{k}\frac{\sqrt{x}}{x-1}  \label{R_as_sum_Gegenbauers} \\
&&\times \sum_{s=1}^{k}\left( -\frac{x+1}{x-1}\right) ^{-s}\binom{k}{s}%
P_{s-1}^{\left( 1,1\right) }\left( \frac{x+1}{x-1}\right) .  \notag
\end{eqnarray}%
Next, we use the contour integral formula for the Jacobi polynomials: 
\begin{equation*}
P_{n}^{\left( \alpha ,\beta \right) }\left( t\right) =\frac{1}{2\pi i}\oint
\left( 1+\frac{t+1}{2}z\right) ^{n+\alpha }\left( 1+\frac{t-1}{2}z\right)
^{n+\beta }z^{-n-1}dz,
\end{equation*}%
with the integration along a small circle around the zero. (See formula
4.4.1 in \cite{szego67}.) It follows that 
\begin{eqnarray*}
\sum_{s=1}^{k}\left( -t\right) ^{-s}\binom{k}{s}P_{s-1}^{\left( 1,1\right)
}\left( t\right) &=&\frac{1}{2\pi i}\oint \sum_{s=1}^{k}\left( -tz\right)
^{-s}\binom{k}{s}\left( 1+tz+\frac{t^{2}-1}{4}z^{2}\right) ^{s}dz \\
&=&\frac{1}{2\pi i}\oint \left( -tz\right) ^{-k}\left( 1+\frac{t^{2}-1}{4}%
z^{2}\right) ^{k}dz,
\end{eqnarray*}%
where we used the binomial theorem in the last step. This is zero for even $%
k.$ For odd $k$ we calculate: 
\begin{equation*}
\sum_{s=1}^{k}\left( -t\right) ^{-s}\binom{k}{s}P_{s-1}^{\left( 1,1\right)
}\left( t\right) =-\frac{1}{t^{k}}\binom{k}{\frac{k-1}{2}}\left( \frac{%
t^{2}-1}{4}\right) ^{\frac{k-1}{2}}.
\end{equation*}%
Next, we set $t=(x+1)/(x-1)$ in (\ref{R_as_sum_Gegenbauers}). Since 
\begin{equation*}
\left( \frac{1+x}{2\sqrt{x}}\right) ^{k}\frac{\sqrt{x}}{x-1}=\left( \frac{t}{%
\sqrt{t^{2}-1}}\right) ^{k}\frac{\sqrt{t^{2}-1}}{2},
\end{equation*}%
hence, for odd $k,$ 
\begin{equation*}
R_{k}\left( a_{1},a_{2}\right) =-2^{-k}\binom{k}{\frac{k-1}{2}}.
\end{equation*}%
\hfill $\square $

\section{Linear statistics of continuosly differentiable functions}

\label{section_C1}

\subsection{Preliminary remarks}

In this section we study the limit distribution for the centered linear
statistics $\mathcal{N}^{o}\left( f,A_{N}\right) $ and $\mathcal{N}%
^{o}\left( g,B_{N}\right) ,$ when $f$ and $g$ are continuously
differentiable functions. First, consider the case of polynomial $f$ and $g$. Recall that the coefficients $\widehat{f}_{k}$ and $%
\widehat{g}_{k}$ are defined as follows: 
\begin{equation}
\widehat{f}_{k}:=\left\{ 
\begin{array}{cc}
\frac{2}{\pi }\int_{-1}^{1}f\left( x\right) T_{k}\left( x\right) \frac{dx}{%
\sqrt{1-x^{2}}}, & \text{for }k\geq 1, \\ 
\frac{1}{\pi }\int_{-1}^{1}f\left( x\right) \frac{dx}{\sqrt{1-x^{2}}}, & 
\text{for }k=1,%
\end{array}%
\right.
\end{equation}%
and similarly for $\widehat{g}_{k}.$ By the orthogonality of Chebyshev
polynomials, one can write: 
\begin{equation*}
f=\sum_{k=0}^{\infty }\widehat{f}_{k}T_{k} \text{ and }%
g=\sum_{k=0}^{\infty }\widehat{g}_{k}T_{k} ,
\end{equation*}%
and for polynomial $f$ and $g,$ the summations in these series are over a
finite number of terms.

Here is a corollary of Theorem \ref{theorem_lin_statistics_Wigner_general}.

\begin{corollary}
\label{Corollary_CLT_for_polynomials}For the real overlapping Wigner
matrices $A_{N}$ and $B_{N}$, and for polynomial functions $f$ and $g,$ the
random variables $\mathcal{N}^{o}\left( f,A_{N}\right) $ and $\mathcal{N}%
^{o}\left( g,B_{N}\right) $ converge in distribution to a two-variate
Gaussian variable with the covariance 
\begin{equation*}
C\left( f,g\right) =\frac{1}{2}\left[ \frac{d_{2}}{2}\widehat{f}_{1}\widehat{%
g}_{1}\gamma +\left( m_{4}-1\right) \widehat{f}_{2}\widehat{g}_{2}\gamma
^{2}+\sum_{k=3}^{\infty }k\widehat{f}_{k}\widehat{g}_{k}\left( \gamma
\right) ^{k}\right] ,
\end{equation*}%
where $\gamma =\Delta /\sqrt{ab}.$
\end{corollary}

Now let us outline the plan of the proof of Theorem \ref{theorem_lin_statistics_diff_functions}.

(1) Take a sequence of polynomials $P_{f,m}$ and $P_{g,m}$ that approximate $%
f$ and $g,$ respectively, in a suitable norm. Let $\mathcal{W}_{m}$ be the
two-variate Gaussian distribution which is the limit for the joint
distributions of $\mathcal{N}^{o}\left( P_{f,m},A_{N}\right) $ and $\mathcal{%
N}^{o}\left( P_{g,m},B_{N}\right) $ as $N\rightarrow \infty .$ Show that the
sequence $\mathcal{W}_{m}$ converges to a limit, a Gaussian distribution $%
\mathcal{W},$ as $m\rightarrow \infty .$

(2) Prove that the joint distributions of pairs $\mathcal{N}^{o}\left(
f,A_{N}\right) $ and $\mathcal{N}^{o}\left( g,B_{N}\right) $ form a tight
family with respect to $N.$ Let $\left\{ \mathcal{Y}_{N}\right\} $ denote
this family and let $\mathcal{Y}$ be one of its limit points.

(3) Show that a suitably defined distance between $\mathcal{Y}$ and $%
\mathcal{W}_{m}$ converges to zero as $m\rightarrow \infty .$

From (1) and (3), we can conclude that $\mathcal{Y}$ must coincide with $%
\mathcal{W}.$ Since this is true for every limit point $\mathcal{Y},$ we
will be able to infer that $\mathcal{N}^{o}\left( f,A_{N}\right) $ and $%
\mathcal{N}^{o}\left( g,B_{N}\right) $ converge to $\mathcal{W}$ as $%
N\rightarrow \infty .$

Before proceeding with this plan, let us derive some preliminary
results.

First, we will need some additional facts about expansions in Chebyshev
polynomials. Consider the change of variable $x=\cos \theta ,$ where $\theta
\in \left[ -\pi ,\pi \right] ,$ and define $F\left( \theta \right) =f\left(
\cos \theta \right) $. If $f(x)$ is absolutely continuous on $\left[ -1,1\right]
,$ then $F\left( \theta \right) $ is absolutely continuous on $\left[ -\pi
,\pi \right] $. By a standard property of Chebyshev polynomials,
$T_n\left( \frac{1}{2} (e^{i\theta} + e^{-i\theta}) \right ) = 
\frac{1}{2} (e^{in\theta} + e^{-in\theta})$. Therefore, the coefficients $\widehat{f}_{n}$ in the expansion of $f$
in the series of Chebyshev's polynomials correspond to the Fourier coefficients in the
Fourier expansion of $F\left( \theta \right) $:%
\begin{equation*}
F\left( \theta \right) =\frac{1}{2}\sum_{n=0}^{\infty }\widehat{f}_{n}\left(
e^{in\theta }+e^{-in\theta }\right) .
\end{equation*}

First, we are going to show that if $f$ is continuously differentiable, then 
$\sum_{n=1}^{\infty }n\left\vert \widehat{f}_{n}\right\vert ^{2}<\infty .$
This will show that the entries of the covariance matrix $V,$ defined in the
statement of Theorem \ref{theorem_lin_statistics_diff_functions}, are
finite. In fact, this holds for a more general class of functions,\ namely,
for the continouous embedding of the Sobolev class $W^{1,p}.$

\begin{lemma}
If $f^{\prime }\in L^{p}\left( \left[ -1,1\right] ,dx\right) $ with $p>1,$
then 
\begin{equation*}
\sum_{n=1}^{\infty }n\left\vert \widehat{f}_{n}\right\vert ^{2}\leq
c_{p}\left\Vert f^{\prime }\right\Vert _{p}^2<\infty .
\end{equation*}
\end{lemma}

\textbf{Proof:} If $f^{\prime }\in L^{p}\left( \left[ -1,1\right] ,dx\right) 
$ with $p\geq 1,$ then%
\begin{eqnarray*}
\int \left\vert F^{\prime }\left( \theta \right) \right\vert ^{p}d\theta
&\ll &\int_{-\pi }^{\pi }\left\vert f^{\prime }\left( \cos \theta \right)
\sin \theta \right\vert ^{p}d\theta \\
&\ll &\int_{-1}^{1}\left\vert f^{\prime }\left( x\right) \right\vert
^{p}\left( 1-x^{2}\right) ^{p/2}\frac{dx}{\sqrt{1-x^{2}}} \\
&\ll &\int_{-1}^{1}\left\vert f^{\prime }\left( x\right) \right\vert ^{p}dx.
\end{eqnarray*}%
(The notation $A(f) \ll B(f)$ for two non-negative functionals $A$ and $B$ means that there is a constant $c$, independent of $f$, such that
$A(f) \leq c B(f)$.)

Hence, $F^{\prime }\in L^{p}\left( \left[ -\pi ,\pi \right] ,d\theta \right)
.$ Moreover, since the interval is finite, hence $F^{\prime }\in L^{s}$ if $%
1\leq s\leq p.$

Recall that the Fourier coefficients of $F^{\prime }$ are $\frac{i}{2}n%
\widehat{f}_{n}.$ Take an $s\in \left( 1,\min (2,p)\right) $ and define $%
r:=s/\left( 2-s\right) >1.$ Then by the H\"{o}lder inequality, 
\begin{equation*}
\sum_{n=1}^{\infty }n\left\vert \widehat{f}_{n}\right\vert ^{2}\leq \left(
\sum_{n=1}^{\infty }\frac{1}{n^{r}}\right) ^{1/r}\left( \sum_{n=1}^{\infty
}\left\vert n\widehat{f}_{n}\right\vert ^{2q}\right) ^{1/q},
\end{equation*}%
where $q=r/\left( r-1\right) =\frac{1}{2}s/\left( s-1\right) .$ The first series on the r.h.s of this inequality is 
convergent because $r>1$. Since $2q>1,$
 the Hausdorff-Young inequality is applicable, and 
\begin{equation*}
\left( \sum_{n=1}^{\infty }\left\vert n\widehat{f}_{n}\right\vert
^{2q}\right) ^{1/\left( 2q\right) }\leq \left( \frac{1}{2\pi }\int
\left\vert F^{\prime }( x) \right\vert ^{s}dx\right) ^{1/s}\ll
\left\Vert f^{\prime }\right\Vert _{s}.
\end{equation*}%
It follows that 
\begin{equation}
\sum_{n=1}^{\infty }n\left\vert \widehat{f}_{n}\right\vert ^{2}\leq
c\left\Vert f^{\prime }\right\Vert _{s}^{2}\leq \widetilde{c}\left\Vert
f^{\prime }\right\Vert _{p}^{2}<\infty .  \label{bound_K_norm}
\end{equation}%
\hfill $\square $

\begin{lemma}
For sequences $x:=\left\{ x_{k}\right\} _{k=1}^{\infty }$ and $y=\left\{
y_{k}\right\} _{k=1}^{\infty },$ define $\left\langle x,y\right\rangle _{\kappa}:=%
\frac{d_{2}}{2}x_{1}\overline{y}_{1}+\left( m_{4}-1\right) x_{2}\overline{y}%
_{2}+\sum_{n=3}^{\infty }nx_{n}\overline{y}_{n}$ and $\left\Vert
x\right\Vert _{\kappa}^{2}:=\left\langle x,x\right\rangle .$ Then $\left\Vert
\cdot \right\Vert _{\kappa}$ is a Hilbert norm induced by the scalar product $%
\left\langle x,y\right\rangle _{\kappa}.$
\end{lemma}

Proof is by verification that $\left\langle c,d\right\rangle _{\kappa}$ is a
scalar product.

Recall that $\mathcal{F}$ is the class of functions continuously differentiable on the interval 
$I_{\delta}=\left[ -1-\delta ,1+\delta \right],$ 
that grow no faster than a polynomial at infinity. For functions $f\in \mathcal{F}$, we define $\left\Vert f\right\Vert _{\kappa}:=\left\Vert
\left\{ \widehat{f}_{i}\right\} _{i=1}^{\infty }\right\Vert _{\kappa},$ where $%
\widehat{f}_{i}$ is the coefficients of the expansion of $f$ in Chebyshev's
polynomials. This is a seminorm on $\mathcal{F}.$ (It is zero on the
subspace spanned by constants$.$)

\subsection{Proof of Theorem \ref{theorem_lin_statistics_diff_functions}}
First, let us approximate the derivative $f^{\prime }$ by polynomials $\widetilde{P}%
_{f,m}\left( x\right) $ of degree $m$. So we take 
$\widetilde{P}_{f,m}\left( x\right) $ so that 

\begin{equation}
\sup_{x\in I_{\delta }}\left\vert f^{\prime }\left( x\right) -\widetilde{P}%
_{f,m}\left( x\right) \right\vert \leq \varepsilon _{m}\rightarrow 0\text{
as }m\rightarrow \infty .  \label{L2_approximation_derivative}
\end{equation}%
Define 
\begin{equation*}
P_{f,m}\left( x\right) :=f\left( -1-\delta \right) +\int_{-1-\delta }^{x}%
\widetilde{P}_{f,m}( t) dt.
\end{equation*}%
Then, we have 
\begin{equation*}
\left\vert f\left( x\right) -P_{f,m}\left( x\right) \right\vert \leq
\int_{-1-\delta }^{x}\left\vert f^{\prime }\left( t\right) -\widetilde{P}_{f,m}( t) \right\vert dt\leq c\varepsilon _{m},
\end{equation*}%
where $c$ is a constant. Hence 
\begin{equation}
\sup_{x\in I_{\delta }}\left\vert f\left( x\right) -P_{f,m}\left( x\right)
\right\vert \leq c\varepsilon _{m}\rightarrow 0\text{ as }m\rightarrow
\infty .  \label{bound_uniform_norm}
\end{equation}%
From (\ref{L2_approximation_derivative}) and (\ref{bound_uniform_norm}) it
follows that $\left\Vert f-P_{f,m}\right\Vert _{Lip}\rightarrow 0,$ where $%
\left\Vert \cdot \right\Vert _{Lip}$ is the Lipschitz norm on the interval $%
I_{\delta }$. (For differentiable functions, Lipschitz norm is defined by $%
\left\Vert h\right\Vert _{Lip}:=\sup_{x\in I_{\delta }}\left\vert h\left(
x\right) \right\vert +\sup_{x\in I_{\delta }}\left\vert h^{\prime }\left(
x\right) \right\vert .$) 

In addition, $\|f'-P_{f,m}^{\prime}\|_p$ is bounded by $\|f'-P_{f,m}^{\prime}\|_{\infty}$. This implies, by (\ref{bound_K_norm}) and (\ref%
{L2_approximation_derivative}), that 
\begin{equation*}
\left\Vert f-P_{f,m}\right\Vert _{\kappa}<c\varepsilon _{m}\rightarrow 0\text{ as 
}m\rightarrow \infty .
\end{equation*}

In particular, by the triangle inequality, $\left\vert \left\Vert
f\right\Vert _{\kappa}-\left\Vert P_{f,m}\right\Vert _{\kappa}\right\vert \rightarrow
0 $ as $m\rightarrow \infty ,$ and also (since $\left\Vert f\right\Vert
_{\kappa}+\left\Vert P_{f,m}\right\Vert _{\kappa}$ $\leq 3\left\Vert f\right\Vert _{\kappa}$
for sufficiently large $m$), $\left\vert \left\Vert f\right\Vert
_{\kappa}^{2}-\left\Vert P_{f,m}\right\Vert _{\kappa}^{2}\right\vert \rightarrow 0$ as 
$m\rightarrow \infty .$

We define $P_{g,m}$ similarly.

By Corollary \ref{Corollary_CLT_for_polynomials}, as $N\rightarrow \infty ,$ 
$\mathrm{Tr}P_{f,m}\left( \widetilde{A}_{N}\right) $ and $\mathrm{Tr}%
P_{g,m}\left( \widetilde{B}_{N}\right) $ converge in distribution to a
bivariate Gaussian variable $\mathcal{W}_{m}$ with the covariance matrix $%
V_{m},$ where%
\begin{eqnarray*}
\left( V_{m}\right) _{11} &=&\frac{1}{2}\left\Vert P_{f,m}\right\Vert
_{\kappa}^{2} \\
\left( V_{m}\right) _{12} &=&\frac{1}{2}\left( \frac{d_{2}}{2}\widehat{%
\left( P_{f,m}\right) }_{1}\widehat{\left( P_{g,m}\right) }_{1}\gamma
+\left( m_{4}-1\right) \widehat{\left( P_{f,m}\right) }_{2}\widehat{\left(
P_{g,m}\right) }_{2}\gamma ^{2}+\sum_{k=3}^{m}k\widehat{\left(
P_{f,m}\right) }_{k}\widehat{\left( P_{g,m}\right) }_{k}\gamma ^{k}\right) \\
\left( V_{m}\right) _{22} &=&\frac{1}{2}\left\Vert P_{g,m}\right\Vert
_{\kappa}^{2}.
\end{eqnarray*}

The diagonal entries $\left( V_{m}\right) _{11}$ and $\left( V_{m}\right) _{22}$ converge to 
$\frac{1}{2}\| f\|_{\kappa}^2$ and $\frac{1}{2}\| g\|_{\kappa}^2$, which are
the diagonal entries of the
matrix $V,$ defined in the statement of Theorem \ref%
{theorem_lin_statistics_diff_functions}. The off-diagonal term $\left(
V_{m}\right) _{12}$ can be written as 
\begin{equation*}
\left\langle P_{f,m}-f,P_{g,m}^{\left( \gamma \right) }\right\rangle
_{\kappa}+\left\langle f,P_{g,m}^{\left( \gamma \right) }-g^{\left( \gamma
\right) }\right\rangle _{\kappa}+\left\langle f,g^{\left( \gamma \right)
}\right\rangle _{\kappa},
\end{equation*}%
where 
\begin{eqnarray*}
P_{g,m}^{\left( \gamma \right) }\left( x\right) &:&=\sum_{k=1}^{m}\widehat{%
\left( P_{g,m}\right) }_{k}\gamma ^{k}T_{k}\left( x\right) ,\text{ and } \\
g^{\left( \gamma \right) }\left( x\right) &:&=\sum_{k=1}^{m}\widehat{g}%
_{k}\gamma ^{k}T_{k}\left( x\right) .
\end{eqnarray*}%
The first two terms are small by the application of the Schwartz inequality
for the scalar product $\left\langle \cdot ,\cdot \right\rangle _{\kappa},$ and
the third term coincides with $V_{12}.$ Hence we can conclude that $%
\left\Vert V_{m}-V\right\Vert $ converges to zero as $m\rightarrow \infty .$
This implies that the Gaussian distributions $\mathcal{W}_{m}$ converge to
the Gaussian distribution $\mathcal{W}$ with the covariance matrix $V.$ This
finishes the first step of the proof.

In order to prove tightness for the family of joint distributions of $\mathcal{N}^{o}\left( f,A_{N}\right)$
  and $\mathcal{N}^{o}\left( g,B_{N}\right) $
(with respect to parameter $N$), we are going to prove that the norms of their
covariance matrices are bounded. In fact, it is enough to prove that
variances of each of $\mathcal{N}^{o}\left( f,A_{N}\right) $ and $\mathcal{N}%
^{o}\left( g,B_{N}\right) $ are bounded, since then the covariance will be
bounded automatically.

Here, we rely heavily on the Poincare inequality property (\textquotedblleft
PI\textquotedblright ) of the matrix entries. The essential feature of the
PI property is that it is well behaved with respect to taking the product of
measures. By definition, the measure $\eta $ on $\mathbb{R}$ has the PI
property, if for some $c_{\eta }>0$ and all differentiable functions $f,$%
\begin{equation*}
\mathbb{V}\mathrm{ar}_{\eta }\left( f\right) \leq c_{\eta }\int \left\vert
f^{\prime }\left( x\right) \right\vert \eta \left( dx\right) .
\end{equation*}%
Then, if $\eta _{K}=\otimes _{i=1}^{K}\eta _{i}$ with $\eta _{i}=\eta $ and
if $h:\mathbb{R}^{K}\rightarrow \mathbb{R}$ is a differentiable function,
then 
\begin{equation*}
\Var_{\eta _K}\left( h\right) \leq c_{\eta }\int
\left\Vert \nabla h\left( x\right) \right\Vert \eta_K \left( dx\right) .
\end{equation*}%
By approximation, this can be further extended to the case when $h$ is
Lipschitz. In particular, if $h$ is a Lipshitz function on $\R^K$,
then we have 
\begin{equation}
\mathbb{V}\mathrm{ar}_{\eta _{K}}\left( h\right) \leq c_{\eta }\left\Vert
h\right\Vert _{Lip}.  \label{Poincare_inequality_2}
\end{equation}

Next, recall that $\mathcal{N}\left( f,A_{N}\right) =\mathrm{Tr}f\left(
A_{N}/\sqrt{4a^{\left( N\right) }}\right) .$ By using the facts about the behavior of the
PI property with respect to scaling and taking products, we find that the joint
distribution of the matrix entries of $A_{N}/\sqrt{4a^{\left( N\right) }}$ 
satisfies the PI property with the constant $c/a^{\left( N\right) }.$ At the same
time, if the function $f\left( x\right) $ is Lipschitz on $\mathbb{R}$, then
the function $\mathrm{Tr}f\left( X\right) $ is Lipshitz on the space of $a_N$-by-$a_N$ 
Hermitian matrices, and 
\begin{equation*}
\left\Vert \mathrm{Tr}f\right\Vert _{Lip}\leq c\sqrt{a_N}\left\Vert
f\right\Vert _{Lip}.
\end{equation*}%
(See Lemma 1.2. in \cite{guionnet_zeitouni00}). Hence, by using (\ref%
{Poincare_inequality_2}), we find that 
\begin{equation}
\mathbb{V}\mathrm{ar}\left( \mathrm{Tr}f\left( A_{N}/\sqrt{4a^{\left(
N\right) }}\right) \right) \leq Cc_{\eta }\left\Vert f\right\Vert _{Lip},
\label{variance_bound_Lipschitz}
\end{equation}%
where $C$ is an absolute constant and $c_{\eta }$ depends only on the
distribution of matrix entries.

A complication arises since under our assumptions, $f$ is assumed Lipschitz
only on the interval $I_{\delta }=\left[ -1-\delta ,1+\delta \right] $.
Outside of $I_{\delta },$ we only know that it has a polynomial growth.
In order to handle this complication, we can write $f$ as a sum of two functions: $%
f=f_{1}+f_{2}$, with $f_{1}$ Lipschitz and bounded everywhere on $\mathbb{R}$%
, $\left\Vert f_{1}\right\Vert _{Lip}<\infty ,$ and $f_{2}$ vanishing on $%
I_{\delta /2}:\left[ -1-\delta /2,1+\delta /2\right] $ and having a
polynomial growth. Then (\ref{variance_bound_Lipschitz}) can be applied to
bound $\mathbb{V}\mathrm{ar}\left( \mathrm{Tr}f_{1}\left( A_{N}/\sqrt{%
4a^{\left( N\right) }}\right) \right) .$

In addition, from the results about the spectra of Wigner matrices, it is
known that the probability for $A_{N}/\sqrt{4a^{\left( N\right) }}$ to have
an eigenvalue outside of $I_{\delta /2}$ becomes exponentially small in $N,$
as $N$ grows. This implies that 
\begin{equation*}
\mathbb{E}\left[ \mathrm{Tr}f_{2}\left( A_{N}/\sqrt{4a^{\left( N\right) }}%
\right) \right] ^{2}\rightarrow 0,\text{ as }N\rightarrow \infty .
\end{equation*}

Since for two random variables, $\xi _{1}$ and $\xi _{2},$ it is true that $%
\sqrt{\mathbb{V}\mathrm{ar}\left( \xi _{1}+\xi _{2}\right) }\leq \sqrt{%
\mathbb{V}\mathrm{ar}\left( \xi _{1}\right) }+\sqrt{\mathbb{V}\mathrm{ar}%
\left( \xi _{2}\right) },$ we can conclude that 
\begin{equation*}
\underset{N\rightarrow \infty }{\lim \sup }\mathbb{V}\mathrm{ar}\left( 
\mathrm{Tr}f\left( A_{N}/\sqrt{4a^{\left( N\right) }}\right) \right) \leq
c\left\Vert f\right\Vert _{Lip},
\end{equation*}%
where the Lipschitz norm is taken over the interval $I_{\delta }.$

A similar argument holds for the random variable $\mathrm{Tr}g\left( B_{N}/%
\sqrt{4b^{\left( N\right) }}\right) ,$ and therefore the norm of the
covariance matrices of these two random variables is bounded. This shows
that the joint distributions of the pairs $\mathcal{N}^{o}\left(
f,A_{N}\right)$ and $\mathcal{N}^{o}\left( g,B_{N}\right) $ form a tight
family and concludes the second step of the proof.

Next, let $\mathcal{Y}$ be a limit point for the distributions $\mathcal{Y}%
_{N}$ of $\{ \mathcal{N}^{o}\left( f,A_{N}\right),$  $\mathcal{N}^{o}\left( g,B_{N}\right) \} ,$
 so that $\mathcal{Y}%
_{N_{k}}\rightarrow \mathcal{Y}$ in distribution for a sequence of $N_{k}$.
We are going to estimate the difference between the characteristic functions
of the distributions $\mathcal{Y}$ and $\mathcal{W}_{m}.$

For convenience, we will assume that all relevant random variables are
realized on a single probability space so that convergence in distribution
reflects convergence almost surely. In this realization, let $Y_{N_{k}}$ and 
$W_{m,N_{k}}$ denote (two-dimensional) random variables that have the same joint distribution
as $\left\{ \mathcal{N}^{o}\left( f,A_{N_{k}}\right) ,\mathcal{N}^{o}\left(
g,B_{N_{k}}\right) \right\} $ and $\left\{ \mathcal{N}^{o}\left(
P_{f,m},A_{N_{k}}\right) ,\mathcal{N}^{o}\left( P_{g,m},B_{N_{k}}\right)
\right\} .$ The variables $Y_{N_{k}}$ \ and $W_{m,N_{k}}$ converge almost
surely to random variables $Y$ and $W_{m},~$that have the distributions $%
\mathcal{Y}$ and $\mathcal{W}_{m}$, respectively. Let $t=(t_1,t_2)\in \R$. Then, 
\begin{eqnarray*}
\left\vert \mathbb{E}e^{itY}-\mathbb{E}e^{itW_{m}}\right\vert &=&\left\vert 
\mathbb{E}\exp \left( it\lim_{N_{k}\rightarrow \infty }Y_{N_{k}}\right) -%
\mathbb{E}\exp \left( it\lim_{N_{k}\rightarrow \infty }W_{m,N_{k}}\right)
\right\vert \\
&\leq &\underset{N\rightarrow \infty }{\lim \sup }\left\vert \mathbb{E}\exp
\left( itY_{N_{k}}\right) -\mathbb{E}\exp \left( itW_{m,N_{k}}\right)
\right\vert \\
&=&\underset{N\rightarrow \infty }{\lim \sup }\left\vert \mathbb{E}\exp
\left( it\left( Y_{N_{k}}-W_{m,N_{k}}\right) \right) -1\right\vert ,
\end{eqnarray*}%
where the inequality follows from Fatou's lemma.

By using (\ref{variance_bound_Lipschitz}), we have 
\begin{equation*}
\underset{N_{k}\rightarrow \infty }{\lim \sup }\mathbb{V}\mathrm{ar}\left[ 
\mathrm{Tr}f\left( A_{N}/\sqrt{4a^{\left( N\right) }}\right) -\mathrm{Tr}%
P_{f,m}\left( A_{N}/\sqrt{4a^{\left( N\right) }}\right) \right] \leq
c\left\Vert f-P_{f,m}\right\Vert _{Lip},
\end{equation*}%
which implies that for the first component of the vector $%
Y_{N_{k}}-W_{m,N_{k}}$ we have the following bound: 
\begin{equation}
\underset{N_{k}\rightarrow \infty }{\lim \sup }
\E\left[ (Y_{N_{k}}-W_{m,N_{k}})_1\right] ^{2}\leq c\left\Vert f-P_{f,m}\right\Vert
_{Lip}.  \label{dif_f_Pf}
\end{equation}%
A similar expression can be written for the second component of $%
Y_{N_{k}}-W_{m,N_{k}}$.

Now, let $\xi _{1}:=(Y_{N_{k}}-W_{m,N_{k}})_1$ and $\xi _{2}:=(Y_{N_{k}}-W_{m,N_{k}})_2$. Note that $\mathbb{E}\xi _{1}=\mathbb{E}\xi _{2}=0.$ Then,

\begin{eqnarray*}
\left\vert \mathbb{E}\exp \left( is\left( t_{1}\xi _{1}+t_{2}\xi _{2}\right)
\right) -1\right\vert &\leq &2s^{2}\mathbb{V}\mathrm{ar}\left( t_{1}\xi
_{1}+t_{2}\xi _{2}\right) \\
&\leq &2s^{2}\left( t_{1}\sqrt{\mathbb{V}\mathrm{ar}\xi _{1}}+t_{2}\sqrt{%
\mathbb{V}\mathrm{ar}\xi _{2}}\right) ^{2},
\end{eqnarray*}%
where the first inequality is a consequence of inequality II.3.14 on p.278
in Shiryaev \cite{shiryaev96}. By using (\ref%
{dif_f_Pf}) and its analogue for the function $g,$ we find: 
\begin{equation*}
\underset{N_{k}\rightarrow \infty }{\lim \sup }\left\vert \mathbb{E}\exp %
\left[ it\left( Y_{N_{k}}-W_{m,N_{k}}\right) \right] -1\right\vert \leq
c\left\Vert t\right\Vert ^{2}\max \left\{ \left\Vert f-P_{f,m}\right\Vert
_{Lip},\left\Vert g-P_{g,m}\right\Vert _{Lip}\right\} ,
\end{equation*}%
which implies that 
\begin{equation*}
\left\vert \mathbb{E}e^{itY}-\mathbb{E}e^{itW_{m}}\right\vert \leq
c\left\Vert t\right\Vert ^{2}\max \left\{ \left\Vert f-P_{f,m}\right\Vert
_{Lip},\left\Vert g-P_{g,m}\right\Vert _{Lip}\right\} .
\end{equation*}%
By our choice, as $m\rightarrow \infty ,$ $P_{f,m}$ and $P_{g,m}$ converge
to $f$ and $g$, respectively, in the Lipschitz norm. Hence, the random
variables $W_{m}$ converge in distribution to $Y.$ This concludes the third
and final step of the proof. As explained before, these three steps imply
that $\left\{ \mathcal{N}^{o}\left( f,A_{N}\right) ,\mathcal{N}^{o}\left(
g,B_{N}\right) \right\} $ converge in distribution to the Gaussian random
variable $W.$ $\square $

\section{Conclusion}

\label{section_conclusion}

We computed the joint distribution of the eigenvalue statistics for two
models of overlapping random matrices. For both the Wigner and sample covariance cases, we
found that the covariance matrix for linear statistics of Chebyshev's $T$-polynomials has
the diagonal structure, and that its diagonal entries depend polynomially on the matrix overlap. 

The computed covariances are different from those found in Borodin's paper for Gaussian matrices. However, the covariances of linear statistics of Chebyshev's $T$-polynomials are the same as in case of the Gaussian matrices provided that the degree of the polynomials is higher than 2 in the Wigner case and higher than 1 in the sample covariance case.

For matrices whose entries satisfy the Poincare inequality property, we extended the results to all continuously differentiable functions.

\appendix

\section{\textsc{Proof of Theorem \ref{theorem_NBT_T_Wishart}:}}
 Define 
\begin{equation*}
S_k :=\sum A_2(\gamma),
\end{equation*}%
where the sum over all closed non-backtracking tailless (``NBT'') paths  $\gamma$ of length $2k$
that start from a vertex in $V.$ Define $\widetilde{S}_k$ similarly except that every path in the sum must start with a vertex in $W$. 

By Theorem  \ref{theorem_NBT_U_Wishart},
\begin{equation*}
\Tr[F_k(AA^{\ast})]=\sum A_2(\gamma),
\end{equation*}
where the sum over all closed non-backtracking paths of length $2k$
that start from a vertex in $V.$ We partition the sum depending on the tail length to get,
\begin{equation*}
\Tr[F_k(AA^{\ast})]=\sum A_2(\gamma)+ \sum A_2(\gamma) + \ldots , 
\end{equation*}
where the first sum on the r.h.s. is the sum over all closed non-backtracking paths with a tail of length 0, the second term is the sum overall non-backtracking paths with a tail of length 1, etc. 

The first term on the r.h.s. is $S_k$. The second term is $(d-1)\widetilde{S}_{k-1}$. Indeed, a tail always contributes the factor of 1 to the product, hence the second term equals the sum over all NBT paths of length $2(k-1)$ that start from a vertex in $W$ multiplied by the number of valid choices for the tail. By a similar counting, we find that the third term on the r.h.s. equals $(c-1)d S_{k-2}$, the fourth term equals $(d-1)cd\widetilde{S}_{k-3}$, the fifth term equals $(c-1)cd^2 S_{k-4}$, etc. Therefore, 
\begin{eqnarray*}
\Tr[F_k(AA^{\ast})]&=&S_k +(d-1)\widetilde{S}_{k-1}+(c-1)dS_{k-2} \\
&&+( d-1) cd\widetilde{S}_{k-3}+( c-1) cd^2 S_{k-4} +...
\end{eqnarray*}

Note that the shift transformation 
\begin{equation*}
\gamma=(v_0,w_{1},v_{1},\ldots
,w_{k},v_0) \to \gamma^{\prime}=( w_{1},v_{1},\ldots
,w_{k},v_{0},w_{1}) 
\end{equation*} defines a bijection of the NBT paths that start with a vertex in $V$ to NBT
paths that start with a vertex in $W,$ and
\begin{equation*}
A_2(\gamma)=A_2(\gamma^{\prime}).
\end{equation*}
Hence, $\widetilde{S}_k=S_k $ for all $k,$ and

\begin{eqnarray*}
\Tr[F_k(AA^{\ast})] &=&\left[ S_k +cd S_{k-2} +(cd)^2 S_{k-4} +\ldots \right] \\
&&+( d-1) \left[S_{k-1}+cd S_{k-3} +\ldots \right] \\
&&-d\left[ S_{k-2}+cd S_{k-4} +\ldots \right] .
\end{eqnarray*}%
In order to infer the dependence of $S_k$ on $A$ we define the polynomial $S_{k}\left( x\right) $
by the formula 
\begin{equation*}
S_{k}\left( x\right) =2\widetilde{T}_{k}\left( x\right) +\frac{\left(
c-d\right) \left( -d\right) ^{k}+cd-1}{d+1}
\end{equation*}%
for $k\geq 1$ and $S_{k}\left( x\right) =0$ for $k\leq 0.$ Then, by using equations (\ref{relation_U_and_T}) and (\ref{relation_F_and_U}), we check that: 
\begin{eqnarray*}
F_k(x) &=&\left[ S_{k}\left( x\right) +cdS_{k-2}\left(
x\right) +\left( cd\right) ^{2}S_{k-4}\left( x\right) +\ldots \right] \\
&&+\left( d-1\right) \left[ S_{k-1}\left( x\right) +cdS_{k-3}\left( x\right)
+\ldots \right] \\
&&-d\left[ S_{k-2}\left( x\right) +cdS_{k-4}\left( x\right) +\ldots \right] .
\end{eqnarray*}%
  This implies the statement of the theorem. \hfill $\square $

\bibliographystyle{plain}
\bibliography{comtest}

\begin{thebibliography}{10}

\bibitem{anderson_guionnet_zeitouni10}
Greg~W. Anderson, Alice Guionnet, and Ofer Zeitouni.
\newblock {\em An Introduction to Random Matrices}, volume 118 of {\em
  Cambridge studies in advanced mathematics}.
\newblock Cambridge University Press, 2009.

\bibitem{anderson_zeitouni06}
Greg~W. Anderson and Ofer Zeitouni.
\newblock A \mbox{CLT} for a band matrix model.
\newblock {\em Probability Theory and Related Fields}, 134:283--338, 2006.

\bibitem{baryshnikov01}
Yuliy Baryshnikov.
\newblock Gues and queues.
\newblock {\em Probability Theory and Related Fields}, 119:256--274, 2001.

\bibitem{borodin10}
Alexei Borodin.
\newblock \mbox{CLT} for spectra of submatrices of wigner random matrices.
\newblock \mbox{\href{http://arxiv.org/abs/1010.0898}{arxiv:1010.0898}}, 2009.

\bibitem{cabanal_duvillard01}
Thierry Cabanal-Duvillard.
\newblock Fluctuations de la loi empirique de grandes matrices al\'{e}atoires.
\newblock {\em Ann. Inst. H. Poincare Probab. Statist.}, 3:373--402, 2001.

\bibitem{callan_smiley05}
David Callan and Len Smiley.
\newblock Noncrossing partitions under rotation and reflection.
\newblock \mbox{\href{http://arxiv.org/abs/math/0510447}{arxiv:math/0510447}},
  2005.

\bibitem{costin_lebowitz95}
Ovidiu Costin and Joel~L. Lebowitz.
\newblock Gaussian fluctuations in random matrices.
\newblock {\em Physical Review Letters}, 75:69--72, 1995.

\bibitem{diaconis_evans01}
Persi Diaconis and Steven~N. Evans.
\newblock Linear functionals of eigenvalues of random matrices.
\newblock {\em Transactions of American Mathematical Society},
  353(7):2615--2633, 2001.

\bibitem{diaconis_shahshahani94}
Persi Diaconis and Mehrdad Shahshahani.
\newblock On eigenvalues of random matrices.
\newblock {\em Journal of Applied Probability}, 31:49--62, 1994.

\bibitem{feldheim_sodin10}
Ohad~N. Feldheim and Sasha Sodin.
\newblock A universality result for the smallest eigenvalues of certain sample
  covariance matrices.
\newblock {\em Geometric and Functional Analysis}, 20:88--123, 2010.
\newblock \mbox{\href{http://arxiv.org/abs/0812.1961}{arxiv:0812.1961}}.

\bibitem{forrester_nagao08}
Peter~J. Forrester and Taro Nagao.
\newblock Determinantal correlations for classical projection processes.
\newblock {\em Journal of Statistical Mechanics: Theory and Experiment},
  2011(08):P08011, 2011.
\newblock \mbox{\href{http://arxiv.org/abs/0801.0100}{arxiv:0801.0100}}.

\bibitem{forrester_nordenstam09}
Peter~J. Forrester and Eric Nordenstam.
\newblock The anti-symmetric \mbox{GUE} minor process.
\newblock {\em Moscow Mathematical Journal}, 9:749--774, 2009.

\bibitem{guionnet_zeitouni00}
A.~Guionnet and O.~Zeitouni.
\newblock Concentration of the spectral measure for large matrices.
\newblock {\em Electronic Communications in Probability}, 5:119--136, 2000.

\bibitem{johansson98}
Kurt Johansson.
\newblock On fluctuation of eigenvalues of random \mbox{Hermitian} matrices.
\newblock {\em Duke Mathematical Journal}, 91:151--204, 1998.

\bibitem{johansson_nordenstam06}
Kurt Johansson and Eric Nordenstam.
\newblock Eigenvalues of \mbox{GUE} minors.
\newblock {\em Electronic Journal of Probability}, 11:1342--1371, 2006.

\bibitem{jonsson82}
Dag Jonsson.
\newblock Some limit theorems for the eigenvalues of a sample covariance
  matrix.
\newblock {\em Journal of Multivariate Analysis}, 12:1--38, 1982.

\bibitem{kostov_mf_shapiro09}
Vladimir~P. Kostov, Andrei Martinez-Finkelshtein, and Boris~Z. Shapiro.
\newblock Narayana numbers and \mbox{S}chur-\mbox{S}zego composition.
\newblock {\em Journal of Approximation Theory}, 161:464--476, 2009.
\newblock available at
  \mbox{\href{http://people.su.se/~shapiro/Articles/EVCSS.pdf}{http://people.su.se/~shapiro/Articles/EVCSS.pdf}}.

\bibitem{lytova_pastur09}
A.~Lytova and L.~Pastur.
\newblock Central limit theorems for linear eigenvalue statistics of random
  matrices with independent entries.
\newblock {\em Annals of Probability}, 37:1778--1840, 2009.

\bibitem{metcalfe11}
Anthony~P. Metcalfe.
\newblock Universality properties of \mbox{G}elfand\mbox{-T}setlin patterns.
\newblock {\em Probability Theory and Related Fields}, 155:303--346, 2013.
\newblock \mbox{\href{http://arxiv.org/abs/1105.1272}{arxiv:1105.1272}}.

\bibitem{reed14}
Matthew Reed.
\newblock PhD thesis, University of California, Davis, 2014.

\bibitem{shcherbina11}
M.~Shcherbina.
\newblock Central limit theorems for linear eigenvalue statistics of the
  \mbox{W}igner and sample covariance random matrices.
\newblock {\em Journal of Mathematical Physics, Analysis, Geometry},
  7:176--192, 2011.
\newblock
  \mbox{\href{http://arxiv.org/abs/math-ph/1101.3249}{arxiv:math-ph/1101.3249}}.

\bibitem{shiryaev96}
A.~N. Shiryaev.
\newblock {\em Probability}.
\newblock Springer, second edition, 1996.

\bibitem{soshnikov00}
Alexander Soshnikov.
\newblock The central limit theorem for local linear statistics in classical
  compact groups and related combinatorial identities.
\newblock {\em Annals of Probability}, 28:1353--1370, 2000.

\bibitem{soshnikov00a}
Alexander~B. Soshnikov.
\newblock Gaussian fluctuation for the number of particles in \mbox{A}iry,
  \mbox{B}essel, sine, and other determinantal random point fields.
\newblock {\em Journal of Statistical Physics}, 100:491--522, 2000.

\bibitem{szego67}
Gabor Szeg\mbox{\"o}.
\newblock {\em Orthogonal Polynomials}.
\newblock American Mathematical Society, third edition, 1967.

\bibitem{tao_vu11}
T.~Tao and V.~Vu.
\newblock Random matrices: universality of the local eigenvalue statistics.
\newblock {\em Acta Mathematica}, 206:127--204, 2011.

\bibitem{zee03}
A.~Zee.
\newblock {\em Quantum Field Theory in a Nutshell}.
\newblock Princeton University Press, 2003.

\end{thebibliography}

\end{document}